\documentclass[final]{siamltex}
\usepackage{amsfonts}
\usepackage{graphicx}

\newcommand{\R}{\Bbb{R}}
\newcommand{\C}{{\mathcal{C}}}

\newtheorem{remark}{Remark}[section]
\newtheorem{conjecture}{Conjecture}[section]

\title{A global stability criterion for scalar
functional differential equations}

\author{Eduardo Liz\thanks{Departamento de
Matem\'atica Aplicada II,
E.T.S.I. Telecomunicaci\'on, Universidad de Vigo,
Campus Marcosende, 36280 Vigo, Spain ({\tt eliz@dma.uvigo.es}).}
        \and Victor Tkachenko\thanks{Institute of Mathematics, National
Academy of Sciences of
Ukraine, Tereshchenkivs'ka str. 3, Kiev, Ukraine ({\tt vitk@imath.kiev.ua}).}
\and Sergei Trofimchuk\thanks{Departamento de Matem\'aticas, Facultad de
Ciencias,
Universidad de Chile, Casilla 653, San\-tia\-go, Chile ({\tt
trofimch@uchile.cl}).}}

\begin{document}
\maketitle

\begin{abstract}
\noindent We consider scalar delay differential equations  $x'(t)
= -\delta x(t) + f(t,x_t) \ (*)$ with nonlinear $f$ satisfying a
sort of negative feedback condition combined with a  boundedness
condition.
The well
known Mackey-Glass type equations, equations satisfying the Yorke
condition, and equations with maxima all fall within our
considerations. Here, we establish a criterion for the global
asymptotical stability of a unique steady state to  $(*)$. As an
example, we study Nicholson's blowflies equation, where our
computations support the Smith's conjecture about the equivalence between
global and local asymptotical stabilities in this population model.
\end{abstract}



\pagestyle{myheadings}
\thispagestyle{plain}
\markboth{Liz,  Tkachenko and Trofimchuk}{Global stability criterion}

\section{Introduction} We start by considering the
simple autonomous linear equation
\begin{equation}\label{simple}
x'(t) = - \delta x(t)+ ax(t-h),
\end{equation}
governed by friction ($\delta \geq 0$) and delayed negative
feedback ($a < 0$). Necessary and sufficient conditions  for the
asymptotic stability of (\ref{simple}) are well known \cite{hl}. For example, in the
simplest case $\delta =0$, Eq. (\ref{simple}) is asymptotically stable if and only
if $-ah \in (0,\pi/2)$. If we allow for a variable delay in (\ref{simple}), we
obtain the equation
\begin{equation}\label{nonautonomouslinear}
 x'(t) = - \delta x(t)+ ax(t-h(t)),  \ 0 \leq h(t) \leq h,
\end{equation}
whose stability analysis  is more complicated  than that of the
autonomous case.  Nevertheless several sharp stability conditions
were established for Eq. (\ref{nonautonomouslinear}). The first of them
is due to Myshkis (see  \cite[p. 164]{hl}) and it states that in
the case $\delta =0$ the inequality $- a \sup_{\Bbb R} h(t) < 3/2$
guarantees the asymptotic stability in (\ref{nonautonomouslinear}).
This condition is sharp (this fact was established
by Myshkis himself). In particular, the upper bound $3/2$
can not be increased  to $\pi/2$. Later on,  the result by Myshkis
has been improved by different authors, the most celebrated
extensions are due to Yorke \cite{yorke} and Yoneyama \cite{Yo}
(both for $\delta = 0$).  Finally, the Myshkis condition has been
recently generalized \cite{ilt}  for  $\delta > 0$:  Eq.
(\ref{nonautonomouslinear}) is asymptotically stable if
\begin{equation}\label{global stability condition with d}
-\frac{\delta}{a}\exp{(-h\delta)}> \ln\frac{a^2-
a\delta}{\delta^2+a^2}.
\end{equation}
We note that for every fixed $a,\delta$ and $h >0$ condition
(\ref{global stability condition with d}) is sharp, and in the
limit case $\delta =0$ it coincides with the  Myshkis
condition. Here the sharpness means  that if $a,\delta, h$ do not
satisfy (\ref{global stability condition with d}), then the
asymptotic stability of Eq. (\ref{nonautonomouslinear}) can be
destroyed by an appropriate choice of a periodic delay $h(t)$ (see
\cite[Theorem 4.1]{ilt}). Returning to Eq. (\ref{simple}), we
can observe that (\ref{global stability condition with d})
approximates exceptionally  well the exact stability domain for
(\ref{simple}) given in \cite{hl}: see Fig. \ref{ffig}, where the
domains  of local (dashed line)
and global (solid line) stability are shown in coordinates
$(-a/\delta, \exp(-\delta h))$. When $\delta =0$, we obtain
$3/2$ as an approximation for $\pi/2$.

It is   a rather surprising fact
that the sharp global stability condition (\ref{global stability
condition with d}) works  not only for  linear
equations, but as well for a variety of nonlinear delay differential
equations of the form
\begin {equation} \label{1}
x'(t) = -\delta x(t) +  f(t,x_t), \ \  (x_t(s) \stackrel{def}=
x(t+s), \ s \in [-h,0]),
\end{equation}
where $f: \Bbb R \times \C \to \Bbb R$, $\C\stackrel{def}= C[-h,0],$ is a
measurable functional satisfying the additional condition {\rm \bf(H)} given
below. Due to the rather general form of {\rm \bf(H)}, Eq. (\ref{1})
incorporates, possibly after some transformations, some of the most celebrated delay equations,
such as equations satisfying the Yorke condition \cite{yorke}, equations of Wright
\cite{hl,lprtt}, Lasota-Wazewska, and Mackey-Glass  \cite{COOKE, kuang, mpn}, and equations with
maxima
\cite{ilt, pt}. Solutions to
some of these equations can exhibit chaotic behavior so that the
analysis of their global stability is of great importance --at
least on the first stage of the investigation (see
\cite[p.148]{kuang} for further discussion). As an example, in
Section \ref{Smith} we  consider Nicholson's blowflies
equation, for which our computations support  the conjecture of Smith posed
in \cite{smith}.

Let us explain briefly the nature of our further assumptions. In
part, they are motivated by the sharp stability results for
(\ref{1}) obtained in \cite{yorke} ($\delta =0$) and \cite{ilt}
($\delta > 0$) under the assumption that for some $a<0$ and for
all $\phi \in \C $, the following Yorke condition holds:
\begin{equation}
\label{ycond} a\mathcal{M}(\phi) \leq f(t,\phi)\leq
-a\mathcal{M}(-\phi).
\end{equation}
Here $\mathcal{M}: \C \to \R $ is the monotone continuous functional
(sometimes called the Yorke functional) defined by
$\mathcal{M}(\phi) = \max\{0, \max_{s\in[-h,0]}\phi(s)\}$.  In
general, $f$ satisfying (\ref{ycond}) is nonlinear in $\phi$. On
the other hand, in some sense it has a ``quasi-linear" form (for example,
$f(\phi) = \max_{s\in[-h,0]} \phi(s)$  can be written as
$f(\phi) = \phi(-s_{\phi})$). In particular, $f$ is sub-linear in
$\phi$, which makes impossible the application of the results from
\cite{ilt, yorke} to the strongly nonlinear cases such as the
celebrated Wright equation
\begin {equation} \label{we}
x'(t) = a(1- \exp(-x(t-h))), \ a < 0,
\end{equation}
which is also globally asymptotically stable if $- ah  \in
(0,3/2)$. Roughly speaking, the Yorke $3/2$-stability condition
does not imply the Wright $3/2$-stability result.  Our recent
studies \cite{lprtt} of (\ref{we}) revealed the following
interesting fact: the essential feature of the function $f(x) =
a(1-\exp(-x)))$ in (\ref{we}) allowing
the
extension of the Wright $3/2$-stability result to some other
nonlinearities is the position of the graph of $f$ with respect
to the graph of the rational function $r(x) = ax/(1+bx)$ which
coincides with $f$, $f'$ and $f''$ at $x =0$.
This suggests the idea to consider a
``rational in
$\mathcal{M}$" version of the ``linear in $\mathcal{M}$" condition
(\ref{ycond}) to manage the strongly nonlinear cases of (\ref{1}).
Therefore, we will assume the following conditions {\rm \bf(H)}:
\begin{description}
  \item[{\rm \bf(H1)}] $f: \Bbb R \times \C \to \Bbb R$ satisfies the
Carath\'eodory condition
  (see \cite[p.58]{hl}). Moreover, for every  $q \in \Bbb R$  there exists 
\ $\vartheta(q)\geq 0 $ such that $f(t,\phi) \leq \vartheta(q)$  almost
everywhere on $\Bbb R$ for every $\phi \in C$  satisfying the
inequality $\phi(s) \geq q, s \in [-h,0]$.
  \item[{\rm \bf(H2)}]  There are $b \geq 0, \ a < 0$ such that
\end{description}
\begin{eqnarray}\label{Y}\qquad \ & f(t,\phi) &
\ \geq \ \frac{a\mathcal{M}(\phi)}{1+b\mathcal{M}(\phi)} \ {\rm
for \ all\ } \phi \in
 \C;
 \\ \qquad \
& f(t,\phi) & \  \leq \
\frac{-a\mathcal{M}(-\phi)}{1-b\mathcal{M}(-\phi)}  {\rm \ for \
all\ } \phi \in \C  {\rm \ such \ that\ }\min_{s\in[-h,0]}\phi(s)
> -b^{-1} \in [-\infty,0).\label{Y2}
\end{eqnarray}
{\rm \bf(H)} is a kind of negative feedback condition combined
with a  boundedness condition; they will cause solutions to remain
bounded and to tend to oscillate about zero. Furthermore, {\rm
\bf(H)} implies that $x = 0$ is the unique steady state solution
for Eq. (\ref{1}) with $\delta>0$. On the other hand, {\rm \bf(H)} does
not imply that the initial value problems for (\ref{1}) have a unique
solution. In
any case, the question of uniqueness is not relevant for our purposes.
Notice finally that if {\rm \bf(H2)} holds with
$b=0$ (which is precisely (\ref{ycond})), then {\rm \bf(H1)} is satisfied
automatically with
$\vartheta(q) = -aM(-q)$.

Now we are ready to state the main result of this work:
\begin{theorem}
\label{main} Assume that {\bf (H)} holds and let $x:[\alpha,
\omega) \to \R$ be a solution of (\ref{1}) defined on the maximal
interval of  existence. Then $\omega= + \infty$ and $x$ is bounded
on $[\alpha, +\infty)$. If, additionally, condition (\ref{global
stability condition with d}) holds, then $\lim_{t \to +\infty}
x(t) = 0$. Furthermore, condition (\ref{global stability condition
with d}) is sharp within the class of equations satisfying {\bf
(H)}: for every triple $a < 0, \ \delta
> 0,\  h > 0$  which do not meet (\ref{global
stability condition with d}), there is a nonlinearity $f$ satisfying {\bf (H)}
and  such that the equilibrium $x(t)=0$ of the corresponding Eq.
(\ref{1}) is not asymptotically stable.
\end{theorem}

It should be noticed that in this paper we do not consider the
limit cases when $b=0$ and/or $\delta =0$. When $b=0, \delta >0$,
Theorem \ref{main} was proved in \cite[Theorem 2.9]{ilt}. The
limit case $\delta =0, \ b \geq 0$ can be addressed  by adapting the
proofs  in \cite{lprtt}. Here, due to the elimination of the
friction term $-\delta x$, an additional condition is necessary
(see \cite{ltt} for details). In this latter case, (\ref{global
stability condition with d}) takes the limit form $-ah \in (0,
3/2)$.
\begin{remark} \label{rescaling}
The set of four parameters ($h >0,\delta>0,a < 0,b>0$) can be
reduced. Indeed, the change of time $\tau = \delta t$ transforms
(\ref{1}) into the same form but with $\delta =1$. Finally, since
$\mathcal{M}$ is a positively homogeneous functional
($\mathcal{M}(k\phi) = k\mathcal{M}(\phi)$ for every $k\geq 0$,
$\phi \in \C$), and since the global attractivity property of the
trivial solution of (\ref{1}) is preserved under the simple
scaling $x = b^{-1}y $, the exact value of $b>0$ is not
important and we can assume that $b = 1$. Also, the change of
variables $x = -y$ transforms (\ref{1}) into $y'(t) =-\delta y(t)
+ [-f(t,-y_t)]$ so that it suffices that at least one of the two
functionals $f(t,\phi), -f(t,-\phi)$ satisfy (\ref{Y}) and
(\ref{Y2}).
\end{remark}

To prove Theorem \ref{main}, in Sections \ref{SCHW} and
\ref{pofmr} we will construct and study several one-dimensional
maps which inherit the stability properties of Eq. (\ref{1}). The
form of these maps depends strongly on the parameters: in fact, we
will split the domain of all admissible parameters given by
(\ref{global stability condition with d}) into several disjoint
parts and each one-dimensional map will be associated to a part.
Some of the maps are rather simple and an elementary analysis is
sufficient to study their stability properties. Some~other maps
are more complicated: for example, the proof of Lemma \ref{GSS}
involves the concept of Schwarz derivative, whose definition and
several of its properties are recalled below. Unfortunately,
several important one-dimensional maps appear in an implicit form
and though this form may be simple, its analysis requires
considerable effort.
For the convenience of the reader, the hardest
and most technical parts of our estimations are placed in an
appendix (Section 6). In Section \ref{Smith}, we will show the
significance of the hypotheses {\bf (H)} again by applying Theorem
\ref{main} to the well-known Nicholson blowflies equation.

\newpage

\section{On  the Smith conjecture and
equations with non-positive Schwarzian} \label{Smith}
\subsection{A global stability condition}
\begin{figure}\label{ffig}
\centering
\includegraphics[totalheight=2.3in]{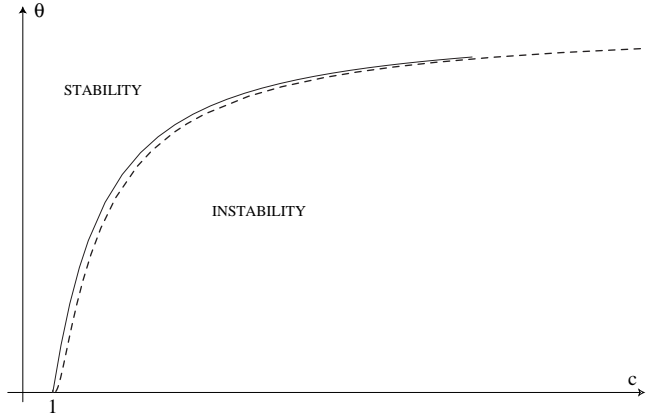}
\caption{Domains of global  and local stability in coordinates
$(c, \theta)$, $c=-a/\delta$, $\theta=\exp(-\delta h)$.}\end{figure}

In this section we will apply our results to the delay
differential equation
\begin{equation} \label{Nich}
N'(u) = -\delta N(u) + p N(u-h) e^{-\gamma N(u-h)}, \ h >  0,
\end{equation}
used by Gurney et al. (see \cite[p.~112]{smith}) to describe the
dynamics of Nicholson's blowflies. Here $p$ is the maximum per
capita daily egg production rate, $1/\gamma$ is the size at which
the population reproduces at its maximum rate, $\delta$ is the per
capita daily adult death rate, $h$ is the generation time and
$N(u)$ is the size of population at time $u.$ In view of the
biological interpretation, we only consider  positive solutions of
(\ref{Nich}). If $p \le \delta,$ Eq. (\ref{Nich}) has only one
constant solution $x \equiv 0.$ For $p > \delta,$ the equation has
an unstable constant solution $x \equiv 0$ and a unique positive
equilibrium $N^* = \gamma^{-1}\ln (p/\delta)$.  Global stability
in Eq. (\ref{Nich}) (when all positive solutions tend to the
equilibrium $N^*$) has been studied by various authors  by using
different methods (see  \cite{COOKE,gyt,smith} for more
references). Nevertheless, the exact global stability condition was not
found. In this aspect, the work \cite{smith}, where the conjecture
about the equivalence between local and global asymptotic
stabilities for Eq. (\ref{Nich}) was posed (see
\cite[p.~116]{smith}), is of special interest for us. Indeed, an
application of our main result to (\ref{Nich}) strongly supports
this conjecture, showing a surprising proximity between  the
boundaries of local and global stability domains; see Fig.
\ref{ffig} and the following proposition:
\begin{theorem}\label{son}
The positive equilibrium $N^*$ of Nicholson's blowflies equation
(\ref{Nich}) is globally asymptotically stable if either $c\in (-1,0]$
or
\begin{equation} \label{globstab}
\mathbf{\theta} > c\ln [(c^2 +c)/(c^2 + 1)]\, , \, c>0,
\end{equation}
where $\mathbf{\theta}=\exp(-\delta h)\, ,\, c = \ln(p/\delta)-1.$
\end{theorem}

It follows from the observation given below (\ref{global stability
condition with d}) that condition (\ref{globstab}) is sharp within
the class of equations $ N'(u) = -\delta N(u) + p N(u-\rho(u))
e^{-\gamma N(u-\rho(u))}$ with variable delay $\rho: \Bbb R \to
[0,h]$.

As it can be seen from (\ref{globstab}), not all parameters
are independent in (\ref{Nich}).  Indeed, if we set $\tau =
h\delta, u = t / \delta, q = p / \delta, x(t)= \gamma N(u)$, then
(\ref{Nich}) takes the form
\begin{equation} \label{LW3}
x'(t) = -x(t) + g(x(t-\tau)),
\end{equation}
where $g(x) = qx\exp(-x)$.
For every $q > 1$, it has a unique positive equilibrium $x(t)
\equiv \ln q$, which is globally asymptotically stable if $\ln(q)\leq 2$ (see \cite{gyt}).
The change of variables $x(t) = \ln q + y(t)$ reduces Eq.
(\ref{LW3}) to the equation $y'(t) = - y(t) + w(y(t-\tau)),$ where
$w(y) = (y + \ln q)e^{-y} - \ln q$. In Section \ref{Pr2}, we will
show that the nonlinearity $w(y)$ satisfies the following
conditions {\bf (W)} within some domain which attracts all
nonnegative solutions of (\ref{LW3}):
\begin{enumerate}
 \item[\bf (W1)] $w \in C^{3}(\R,\R)$, \ $xw(x) < 0$  for  $x\neq 0$ and $w'(0) < 0$.
 \item[\bf (W2)] $w$ is bounded below and has at most one critical point
$x^*\in \R$ which is a local extremum.
 \item[\bf (W3)] The Schwarz derivative $(Sw)(x)=w'''(x)(w'(x))^{-1}-(3/2)
\left(w''(x)(w'(x))^{-1}\right)^2$ of $w$ is nonnegative:
$(Sw)(x)\leq 0$ for all $x\neq x^*$.
\end{enumerate}
Since $w'(0)= \ln(e\delta/p)<0$ and $w''(0)>0$ if $\ln q>2$, Theorem
\ref{son}  is a consequence of the following results:
\begin{lemma}{\rm \cite{lprtt}}\label{lk}  Let $w$ meet
conditions {\bf (W)} and $w''(0) > 0$. Then the functional
$f(t,\phi)= w(\phi(-h))$  satisfies hypotheses {\rm \bf(H)} with
$a = w'(0)$ and $b = -w''(0)/(2w'(0))$.
\end{lemma}
\begin{corollary}
\label{codo} Suppose that $w$ satisfies {\bf (W)} and $w''(0)>0$. If
(\ref{global stability condition with d}) holds with $a =
w'(0)$, then the trivial steady state attracts all solutions of
the delay differential equation
\end{corollary}
\begin{equation}\label{111}
x'(t) = - \delta x(t) + w(x(t-h)),\ \delta > 0.
\end{equation}

Corollary \ref{codo} can be  applied in a similar way to obtain 
 global stability conditions for the positive equilibrium of
other delay differential equations arising in biological models.
For example, we can mention the celebrated Mackey-Glass equation
proposed in 1977 to model blood cell populations (see, e.g.,
\cite{mpn}), which is of the form (\ref{111}) with
$w(x)=b/(1+x^n)$, $b>0$, $n>1$. Another important model that can
be considered within our approach is the Wazewska-Czyzewska and
Lasota equation describing the erythropoietic (red-blood cell)
system. In this case $w(x) = b_1\exp (-b_2 x), \ b_i > 0$.

As it was proved in \cite{lprtt}, the conclusion of Corollary \ref{codo} also holds for $\delta=0$ by
replacing (\ref{global stability condition with d}) with its limit form $-ah\leq 3/2$.
In the
particular case of the Wright equation, this result coincides with the $3/2$-
stability theorem by Wright (see \cite{lprtt} for more details).

\subsection{The Smith and Wright conjectures revisited}
Let us look again on Fig. \ref{ffig}, which shows the boundaries
of the  domains of local and global asymptotic stability for the
Nicholson equation;
this observation (as well as Proposition \ref{sw} stated below) suggests 
the following
\begin{conjecture}
\label{conj} Under conditions {\bf (W)}, the trivial solution of
Eq. (\ref{111}) is globally attracting if it is locally
asymptotically stable.
\end{conjecture}

An interesting particularity of Conjecture \ref{conj} is that it
coincides with the celebrated Wright conjecture if we take $\delta
= 0,\ w(x) = a(1-\exp(-x)),$ and it coincides with  the Smith
conjecture if we take  Nicholson's blowflies equation.

Now, the following result was obtained in
\cite{waln} as a simple consequence of an elegant approach toward
stable periodic orbits for Eq. (\ref{111}) with Lipschitz
nonlinearities:
\begin{proposition}{\rm \cite{waln}} \label{11} For every $\alpha \geq 0$
there exists a smooth strictly decreasing function $w$ satisfying
{\bf (W1)}, {\bf (W2)}, $-w'(0) = \alpha$ and such that Eq.
(\ref{111}) has a nontrivial periodic solution which is
hyperbolic, stable and exponentially attracting with asymptotic
phase (so therefore (\ref{111}) is  not globally stable).
\end{proposition}

Proposition \ref{11} shows clearly that the strong dependence
between local (at zero) and global asymptotical stabilities of Eq.
(\ref{111}) cannot be explained only with the concepts presented
in {\bf (W1), (W2)}. We notice here that the condition of negative
Schwarz derivative in Eq. (\ref{111}) appears naturally also in
some other contexts of the theory of delay differential equations,
see e.g. \cite[Sections 6--9]{mpn}, where it is explicitly used,
and \cite[Theorem 7.2, p.~388]{hl}, where the condition $Sw < 0$
is implicitly required.

\section{Preliminary stability analysis of Eq. (\ref{1})}
\label{SCHW} Throughout  the paper, in view of Remark
\ref{rescaling}, we assume that $\delta =1$ in (\ref{1}) and $b=1$
in (\ref{Y}), (\ref{Y2}). Hence, with $\theta \stackrel{def} = \exp(-h)$,
(\ref{global stability condition with d}), (\ref{1}),  (\ref{Y}) and (\ref{Y2})
take, respectively, the forms
\begin{eqnarray}\label{global stability condition}
 -\theta/a
&>& \ln\frac{a^2-a}{a^2+1};\\
 x'(t)& =& -x(t) +  f(t,x_t); \label{1pr} \\
 f(t,\phi)&\geq& r(\mathcal{M}(\phi)), \ {\rm
for \ all\ } \phi \in
 \C; \label{YY} \\
f(t,\phi)&\leq& r(-\mathcal{M}(-\phi)),   {\rm \ for \ all\ } \phi
\in \C  {\rm \ such \ that\ }\min_{s\in[-h,0]}\phi(s)
> -1, \label{YY2}
\end{eqnarray}
where the rational function $r(x) = ax/(1+x)$ will play a key role
in our constructions.
In this section, we establish that the ``linear" approximation to
(\ref{global stability condition}) of the form
\begin{equation}\label{global stability}
\quad -\theta/a> -(a+1)/(a^2+1)
\end{equation}
implies the global stability of Eq. (\ref{1pr}) (note here that
$\ln(1+x) < x$ is true for $x > 0$).

In the sequel we will use some properties of the Schwarz
derivative. The following lemma can be checked by direct
computation:
\begin{lemma}
\label{compsw} If $g$ and $f$ are functions which are at least
$C^3$ then $S(f\circ g)(x)=(g'(x))^2 (Sf)(g(x))+(Sg)(x)$. As a
consequence, the inverse $f^{-1}$ of a smooth diffeomorphism $f$
with $Sf >0$ has negative Schwarzian: $ Sf^{-1} < 0$.
\end{lemma}

We will also need the following lemma from \cite{dev}:

\begin{lemma}{\rm \cite[Lemma 2.6]{dev}}
\label{lsin} Let   $q: [\alpha,\beta] \to [\alpha,\beta]$ be a
$C^3$ map with $(Sq)(x) < 0$ for all $x$. If $\alpha< \gamma
<\beta$ are consecutive fixed points of some iteration $g=q^N$ of $q$, $N\geq 1,$ and
$[\alpha,\beta]$ contains no critical point of $g$, then
$g'(\gamma)>1$.
\end{lemma}

This lemma allows us to prove the following proposition, which plays 
a central role in
our analysis.

\begin{proposition}
\label{sw}  Let  $q: [\alpha,\beta] \to [\alpha,\beta]$ be a $C^3$
map with a unique fixed point $\gamma$ and with at most one
critical point $x^*$ (maximum).  If $\gamma$ is locally
asymptotically stable and the Schwarzian derivative $(Sq)(x) < 0$
for all $x \not=x^*$, then $\gamma$ is the global attractor of
$q$.
\end{proposition}

\begin{proof}
 Let $W$ be the connected component of the open set   $S= \{x\in  [\alpha,\beta]\, :\,
\lim_{k\to +\infty}q^k(x)= \gamma\}$ which contains $\gamma$. Clearly,
$g(W)\subset W$. If $W\neq  [\alpha,\beta]$, then we have three
possibilities:
 $W=[\alpha,r)$, $W=(l,\beta]$ or $W=(l,r)$, $\alpha <l<r<\beta$.

If  $W = [\alpha,r)$  then $q(r) = \lim_{\epsilon \to 0+} q(r-
\epsilon) \leq r \leq q(r)$, a contradiction with the  fact  that
$q$  does not have fixed points in $ [\alpha,\beta]$ different
from $\gamma$. The case $W = (l,\beta]$ is completely  analogous.

In the case $W = (l,r)$, by the same arguments, it should hold
$q(l)=r$, $q(r)=l$. Thus $l<\gamma <r$ are consecutive fixed
points of $g=q^2$ and $g'(\gamma)=(q'(\gamma))^2\leq 1$. By Lemma
\ref{lsin}, $x^*\in (l,r)$ and therefore $q(x^*)<r$. Since $q$ has
a maximum at $x^*$, $r>q(x^*)\geq q(l)$, a contradiction.

Hence $W= [\alpha,\beta]$ and therefore $\{\gamma\}$ attracts each
point of $ [\alpha,\beta]$. This implies that $\gamma$ is the
global attractor of $q$ (see \cite[Chapter 2]{HMO}).
\end{proof}

Now we are in a position to begin the stability analysis of
Eq. (\ref{1pr}).
\begin{lemma} \label{L1} Suppose that {\rm \bf(H)}  holds and
let $x:[\alpha-h, \omega) \to \R$ be a solution of Eq. (\ref{1pr})
defined on the maximal interval of  existence. Then $\omega= +
\infty$ and $M = \limsup_{t \to \infty}x(t), \ m = \liminf_{t \to
\infty}x(t)$ are finite. Moreover, if $m \geq 0$ or $M \leq 0$,
then $M = m = \lim_{t \to \infty}x(t) =0$.
\end{lemma}
\begin{proof}  Note that (\ref{YY}) implies that $f(t,\phi) \geq a$
for all $t \in \R$ and $\phi \in \C$. Next, if $q \le x_{\alpha}(s)
\le Q,\ s \in [-h, 0]$ then  for all $t \in  [\alpha, \omega)$, we
have
\begin{eqnarray}\label{ozenka snizu}
x(t) &= & \exp(-(t-\alpha))x(\alpha) +
\int^{t}_{\alpha}\exp(-(t-s))f(s,x_s)ds  \\
&\geq& a +(\min \{q, a \}-a)\exp(-(t-\alpha)) \geq \min \{q, a \}.
\nonumber
\end{eqnarray}
Next, {\rm \bf(H1)} implies that $f(s,x_s) \leq \vartheta(\min
\{q, a \})$ for all $s\geq \alpha$, so that
$$x(t) \le \max \{Q , 0 \}+ \vartheta(\min
\{q, a \}), \ t \in [\alpha, \omega).$$ Hence $x(t)$ is bounded on
the maximal interval of existence, which implies the boundedness of the
right hand side of Eq. (\ref{1pr}) along $x(t)$. Thus $\omega =
+\infty$ due to the corresponding continuation theorem (see \cite[Chapter 2]{hl}).

Next, suppose, for example, that $M = \limsup_{t \to
\infty}x(t)\leq 0$. Thus we have $\lim\limits_{t \to
\infty}\mathcal{M}(x_t) =0$ so that,
in virtue of (\ref{YY}),  $f(t,x_t) \geq
\inf\limits_{s \geq t} a\mathcal{M}(x_s)/(1+\mathcal{M}(x_s))
\stackrel{def}= a(t)$ where $a: [\alpha, +\infty) \to (-\infty,0]$
is nondecreasing and continuous, with $\lim_{t \to +\infty} a(t) =
0$. Thus, by (\ref{ozenka snizu}), $x(t) \geq
\exp(-(t-\beta))x(\beta) + a(\beta)$ for all $t \geq \beta >
\alpha$. This implies that $m =
\liminf_{t \to \infty}x(t) =0$ so that $M=0$.
\end{proof}

\begin{lemma}
\label{GT} Suppose that {\rm \bf(H)}  holds and
let $x:[\alpha-h, \infty) \to \R$ be a solution of Eq. (\ref{1pr}). If $x$ has a negative
local minimum at some point $s>\alpha$, then
$\mathcal{M}(x_{s}) > 0$. Analogously, if $x$ has a positive
local maximum at  $t>\alpha$, then
$\mathcal{M}(-x_{t}) > 0$.
\end{lemma}
\begin{proof}
If
$x(u) \leq 0$ for all $u \in [s-h, s]$, then $x'(s) \geq -x(s)
+ r(\mathcal{M}(x_s))> 0$, a contradiction.  The other case is similar.
\end{proof}

\begin{lemma}
\label{GSS}  Suppose that {\rm \bf(H)}  holds and
let $x:[\alpha-h, \infty) \to \R$ be a solution of Eq. (\ref{1pr}). If (\ref{global
stability}) holds then $\lim_{t \to \infty}x(t) =0$.
\end{lemma}
\begin{proof}
Let  $M = \limsup_{t \to \infty}x(t), \ m = \liminf_{t \to
\infty}x(t)$. In view of Lemma \ref{L1}, we only have to consider the case
$m<0<M$, since otherwise ($m\geq 0$ or $M\leq 0$) we have a non
oscillatory solution to Eq. (\ref{1pr}), which tends to zero as $t\to
+\infty$. Thus in the sequel we will only consider the oscillating
solutions $x(t)$. In this case there are two sequences of points
$t_j, s_j$ of local maxima  and local minima respectively such
that $x(t_j)= M_j \to M, x(s_j) = m_j \to m$ and $s_j, t_j \to
+\infty$ as $j \to \infty$.

First we prove that $M=m=0$ if  $a(1-\theta) > -1$. Indeed, for each $s_j$
we can find $\varepsilon_j \to 0+$ such that $0 < \mathcal{M}(x_s)
< M + \varepsilon_j$ for all $s\in [s_j-h,s_j]$. Next, by Lemma
\ref{GT}, there exists $h_j \in [0,h]$ such that $x(s_j-h_j) = 0$.
Therefore, by the variation of constants formula,
$$
m_j = \int_{s_j-h_j}^{s_j}e^{s-s_j}f(s,x_s)ds \geq
\int_{s_j-h_j}^{s_j}e^{s-s_j}r(\mathcal{M}(x_s))ds \geq
r(M+\varepsilon_j)(1-\theta).
$$
As a limit form of this inequality, we get $m \geq
r(M)(1-\theta)$. Hence $m > -1$ and we can use  (\ref{YY2}) for
$\phi = x_t$ with sufficiently large $t$. Thus, in a similar way,
we obtain that
\begin{eqnarray*}
M_j = \int_{t_j-h^*_j}^{t_j}e^{s-t_j}f(s,x_s)ds \leq
\int_{t_j-h^*_j}^{t_j}e^{s-t_j}r(-\mathcal{M}(-x_s))ds  \leq
r(m-\varepsilon^*_j)(1-\theta)
\end{eqnarray*}
for some sequences $\varepsilon^*_j \to 0+$  and $h^*_j \in
[t_j-h,t_j]$. Hence we obtain $M \leq r(m)(1-\theta) \leq
r(r(M)(1-\theta))(1-\theta)$. This gives $M^2 \leq M(a(1-\theta)
-1)$ which is only possible when $M=0$.

Now, assume that
$a(1-\theta) \leq  -1$. Since
$m \geq r(M)(1-\theta) > a$ (see the first part of the proof),
we conclude that $r^{-1}(m) = m/(a-m)
> 0$ is well defined. Next, for $s_j$ we can find a sequence of
positive $\epsilon_j \to 0$ such that $m_j < m +\epsilon_j < 0$.
We claim that $x(s_j-h_j) \geq r^{-1}( m+ \epsilon_j)$ for some
$h_j \in [0,h]$. Indeed, in the opposite case, $x(s) <  r^{-1}( m+
\epsilon_j)$ for all $s$ in some open neighborhood of
$[s_j-h,s_j]$. Thus $f(s,x_s) \geq r(\mathcal{M}(x_s))
> m+ \epsilon_j$ for all $s$ close to $s_j$. Finally,
$x'(s) > - x(s) +  m+ \epsilon_j > 0$ almost everywhere in some
neighborhood of $s_j$, contradicting  the choice of $s_j$.

Next, there exists a sequence of positive $\epsilon^*_j \to 0$
such that $x(s) < M + \epsilon^*_j$ for all $s \in [s_j-h, s_j ].$
 Therefore, by  the variation of constants formula,
\begin{eqnarray*}
m_j &=& x(s_j) = e^{-h_j} x(s_j -h_j) + \int_{s_j-h_j}^{s_j} e^{s-
s_j}f(s, x_s)ds  \\ &\geq& e^{-h_j}r^{-1}( m+ \epsilon_j) +
\int_{s_j-h_j}^{s_j}e^{s-s_j}r(\mathcal{M}(x_s))ds \\ &\geq&
e^{-h_j}r^{-1}( m+ \epsilon_j) +  r(M+\epsilon^*_j)(1-e^{-h_j})
\geq \theta r^{-1}( m+ \epsilon_j) + r(M+\epsilon^*_j)(1-\theta),
\end{eqnarray*}
so that $m - \theta r^{-1}( m)  \geq r(M)(1-\theta) 
\geq a(1-\theta)$. This implies that $m \geq
a(1-\theta)(a-m)/(a-m-\theta)
> -1$, where the last inequality is evident when $1 + a(1-\theta) =0 $ and
follows from the relations $m < 0 \leq (a^2(1-\theta) +a
-\theta)/(1 + a(1-\theta))$ otherwise.
Since $m
> -1$
we can use
(\ref{YY2}) for $\phi = x_t$ with sufficiently large $t$. Thus, in
a similar way, we obtain that
\begin{eqnarray*}
M_j &=& x(t_j) = e^{-h^\#_j} x(t_j -h^\#_j) +
\int_{t_j-h^\#_j}^{t_j} e^{s- t_j}f(s, x_s)ds \leq
e^{-h^\#_j}r^{-1}( M - \varepsilon_j) +\\
&+& \int_{t_j-h^\#_j}^{t_j}e^{s-t_j}r(-\mathcal{M}(-x_s))ds  
 \leq \theta r^{-1}( M-\varepsilon_j) +
r(m-\varepsilon^*_j)(1-\theta)
\end{eqnarray*}
for some sequences $\varepsilon_j, \varepsilon^*_j \to 0+$ and
$h^\#_j \in [t_j-h,t_j]$. 
Thus $\psi(M)\stackrel{def} = M- \theta r^{-1}(M)\leq
r(m)(1-\theta)$. Now, $\psi:(a, +\infty) \to \R$
is a strictly increasing bijection so that $\chi(x)= \psi^{-1}((1
- \theta)r(x))$ is well defined and strictly decreases on $(-1,
+\infty)$. A direct computation shows that $\chi(-1^-)= +\infty$
and that  $\chi(+\infty)= \psi^{-1}((1 - \theta)a) > -1$.
Therefore $\chi: [\chi(+\infty),\chi^2(+\infty)] \to
[\chi(+\infty),\chi^2(+\infty)]$. Moreover, since
$M \le \chi(m), \ m \ge \chi(M),$ we conclude that $m,M \in
[\chi(+\infty),\chi^2(+\infty)]$ and that $[m,M] \subset
\chi([m,M])$. Next,  for $x > a$ we obtain by direct computation that
$(S\psi)(x) = -6\theta a(a^2-2xa+x^2-\theta a)^{-2} >0.$

Since $(Sr)(x)=0$ for all $x>-1$, it follows from  Lemma \ref{compsw}  that
$$(S\chi)(x) = ((1 - \theta)r'(x))^2(S\psi^{-1})((1 -
\theta)r(x)) < 0.$$ Finally, by (\ref{global
stability}), $\chi'(0) =
(1-\theta)a^2/(a-\theta) \in (-1,0)$ so that we  apply
Proposition \ref{sw} (where we set $q = \chi, \ [\alpha,\beta] =
[\chi(+\infty),\chi^2(+\infty)]$ and $\gamma = 0$) to conclude
that $\chi^k( [\alpha,\beta]) \to 0$ as $k \to \infty$.  Since
$[m,M] \subseteq \chi^k([m,M]) \subseteq \chi^k( [\alpha,\beta])$
for all integers $k\geq 1$, it is clear that $m=M=0$.
\end{proof}

\section{Proof of the main result} \label{pofmr}
The analysis done in the previous section shows that the only case that remains
to consider is
when 
$$0 < \ln\frac{a^2-a}{a^2+1} <  -\theta/a \leq -\frac{a+1}{a^2+1}.
$$
This case will be studied in the present section: we start
describing a finer decomposition of the above indicated domain of
parameters (denoted below as $\mathcal{D}$).
\subsection{Notations and domains}
\label{nodo} 
\begin{figure}[ht]
\centering
\includegraphics{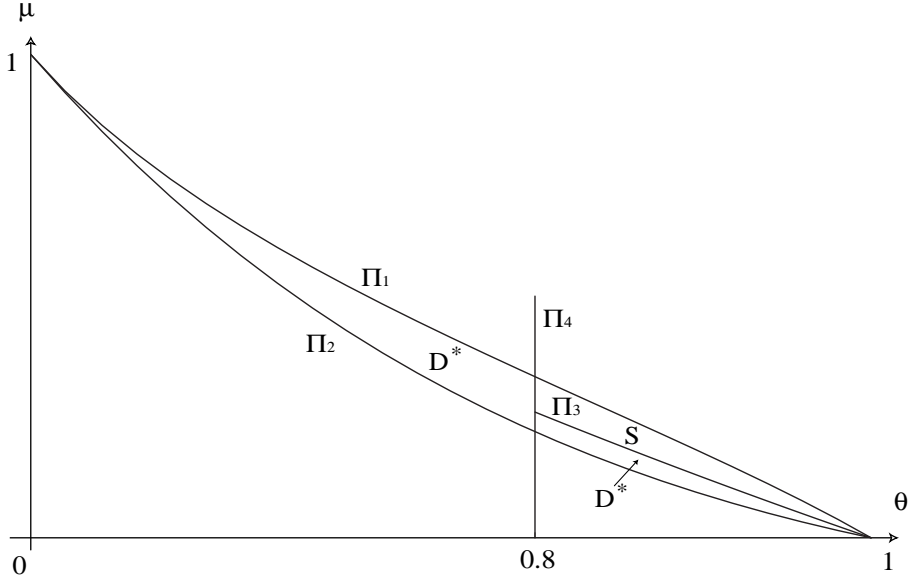}
\caption{Domains of global  stability in coordinates $(\theta,
\mu)$} \label{fig}
\end{figure}

In the sequel, we will always assume  that $h > 0$ and
$a < -1$, and will use the following notations:
\begin{eqnarray*} \label{q}
& & \theta = \exp(- h); \ \ \lambda = \exp(\theta/a); \ \
a_* = a + \frac{\theta}{1-\theta}; \ \ \mu = -\frac1a; \ \
\\[3mm]
& & \alpha(a,\theta)  = (1-a)\exp(\theta/a)+a ; \label{a}\\[3mm]
& & \beta(a,\theta)  =-\frac{a^2 + \exp(\theta/a)(1-2a +2\theta
(a-1))- (1-a)^2\exp(2\theta/a)}{a^2+(a-a^2)\exp(\theta/a)}; \\[3mm]
& & \gamma(a,\theta)  =a^3\alpha(a,\theta)\frac{1-\theta+ \ln
\theta}{2- \theta+\ln \theta}; \ \ \mathcal{R}(r)=
\mathcal{R}(r,a,\theta) = \frac{\alpha(a,\theta) r}{1 -
\beta(a,\theta) r}.
\end{eqnarray*}
Obviously, $\theta, \ \mu 
\in (0,1), \ a_* > a,$ and $\gamma(a,\theta)$ is well defined for
all $\theta\in (0.16,1),$ where it can be checked that $2-
\theta+\ln \theta > 0$. Next, we will need the following four
curves considered within the open square $(\theta, \mu)\in
(0,1)^2$:
$$\begin{array}{l}
\Pi_1 = \{ (\theta, \mu)\, :\, \theta =\Pi_1(\mu) \stackrel{def}=
\displaystyle\frac{1-\mu}{1+\mu^2}\} \; ; \;
\Pi_2 = \{(\theta, \mu)\, :\,\theta = \Pi_2(\mu)\stackrel{def}=
\displaystyle\frac{1}{\mu}\ln\frac{1 +\mu}{1+\mu^2}\} ;\\[3mm]
\Pi_3 = \{(\theta, \mu)\, :\,\theta = \Pi_3(\mu)\stackrel{def} =
\displaystyle\frac{95-108\mu}{5(19+5\mu)}\} \; ; \;
 \Pi_4 = \{(\theta, \mu)\, :\, \theta = \Pi_4(\mu)\stackrel{def}= 0.8\}.
\end{array}
$$
The geometric relations existing between  curves $\Pi_1- \Pi_4$
are shown schematically on Fig. \ref{fig}. 
Notice that all three curves $\Pi_j, \ j \not= 4,$ have the
following asymptotics  at zero:  $\Pi_j(\mu) = 1 -
k_j\mu + o(\mu)$, where  $k_1 = 1, k_2 = 1.5, k_3 =  1.4$.  An
elementary analysis  shows that $\Pi_3$ does not
intersect  $\Pi_1$ and $\Pi_2$  when $\theta \in
(0.8,1)$. Next, to prove our main result, we will have to
use different arguments for the different domains of parameters
$a, h$.  For this purpose, we introduce here the following three
subsets $\mathcal{D}, \mathcal{D}^*, \mathcal{S}$ of $(0,1)^2$
\begin{eqnarray*} \label{qh}
 & & \mathcal{D}= \{(\theta, \mu): \Pi_2(\mu) \leq \theta \leq
\Pi_1(\mu) \} ;\\[3mm]
& & \mathcal{D}^* = \mathcal{D}\setminus \mathcal{S}, \ {\rm
where}\ \mathcal{S}= \{(\theta, \mu) \in \mathcal{D}: \ \theta
\in[0.8, 1), \ \Pi_3(\mu) \leq \theta \leq \Pi_1(\mu)\}.
\end{eqnarray*}
We can see that $\mathcal{D}$ is situated between $\Pi_1$ and
$\Pi_2$, while the sector $\mathcal{S}$ is placed among $\Pi_1,$
 $\Pi_3,$ and $\Pi_4$.
Sometimes it will be more convenient for us to use the coordinates
$(a,\theta)$ instead of $(\theta,\mu)$, we will preserve the same
symbols for the domains and curves considered both in $(a,\theta)$ and
$(\theta,\mu)$.

Let us end this section indicating several useful estimations
which will be of great importance for the proof of our main
result.
\begin{lemma} \label{albetand} We have $\alpha(a,\theta)>0$,
$\beta(a,\theta)> 0$, and $a \alpha(a,\theta)/(1- a\beta(a,\theta)) >
-1 $ for all $(a,\theta) \in \mathcal{D}$. Next, if $(a,\theta)
\in \mathcal{S},$ then $\gamma(a,\theta) < 1.$
\end{lemma}

The proof of the lemma is given in Section 6 (Lemmas
\ref{albet}, \ref{albeta}, \ref{dom}).

\subsection{One-dimensional map $\mathbf{F: I \to }\R$} \label{subF}
Throughout this subsection, we will suppose that $(a,\theta) \in
\mathcal{D}$. Therefore $ a(\theta-1)/\theta -1 > 0$ so that the
interval $I = (-1, a(\theta-1)/\theta -1)$ is not empty.
Furthermore, $t_{1} = t_{1}(z)=-\ln (1-z/r(z)) \in [-h,0]$ for
every $z \in I\setminus\{0\}$.  Consider now the map $F: I \to \R$
defined in the following way:
$$ F(z)= \left\{
\begin{array}{lll}
    0   & {\rm if} \ z = 0; \\
    \min\limits_{t\in [0,h]}y(t,z) & {\rm if} \ z \in I \  {\rm and } \ z > 0; \\
    \max\limits_{t\in [0,h]}y(t,z) & {\rm if} \ z \in (-1,0),
\end{array}\right.
$$
where $y(t,z)$ is the solution of the initial value problem
$y(s,z)=z, \ s \in [t_1(z)-h,t_1(z)], \ z \in I$ for
\begin{equation}
\label{ep} y'(t)=-y(t)+r(y(t-h)).
\end{equation}
Observe that $y(0,z) = 0$ for all $z \in I$ since
$y(t,z)=r(z)(1-\exp{(-t)})$ for all $t \in [t_{1}(z),t_{1}(z)+h]$.
The following lemma explains why we have introduced such $F$
(moreover, condition (\ref{global stability condition})
says precisely that $F'(0)>-1$, see Section 6.2):
\begin{lemma}
\label{44} Let $x(t)$ be a solution of (\ref{1pr}) and set $M =
\limsup_{t \to \infty}x(t),$ $m = \liminf_{t \to \infty}x(t)$. If
$m,M \in I$, then $m \ge F(M)$ and $M \le F(m)$.
\end{lemma}
\begin{proof} Consider two sequences of extremal values $m_j=
x(s_j) \to m$, $ M_j = x(t_j) \to M$ such that  $s_j \to +\infty,
t_j \to +\infty$ as $j \to \infty$. Let $\varepsilon > 0$ be such
that $(m - \varepsilon, M+\varepsilon) \subset I$. Then  $m_j \ge
m-\varepsilon$ and $M_j \le M+\varepsilon$ for big $j$.  We will
prove  that $m \geq F(M)$, the case $M \leq F(m)$ being completely
analogous.

By Lemma \ref{GT},  we can find $\tau_j \in [s_j-h, s_j]$ such
that  $x(\tau_j) = 0$ while $x(t) < 0$ for $t \in (\tau_j,s_j]$.
Next, $v_j = \tau_j+ t_1(M+\varepsilon) \geq \tau_j -h$ because of
$M+\varepsilon \in I$. Thus the solution $y(t)$ of Eq.
(\ref{ep}) with initial condition
$y(s) = M+\varepsilon, \ s \in [v_{j}-h,v_{j}]$,  satisfies $y(\tau_j) = 0$
while $M+\varepsilon \equiv y(t) \geq x(t)$ for all $t \in [v_{j}-h,v_{j}]$.
Furthermore, for all $s \in [v_{j},\tau_j],$ we have $\mathcal{M}(x_s) \le
M+\varepsilon$ so that, by (\ref{YY}), $f(s, x_s) \ge
r(\mathcal{M}(x_s)) \ge r(M+\varepsilon),$ and
$$y(t) - x(t) =  \int_{\tau_j}^{t} e^{-(t-s)}
[r(M+\varepsilon) - f(s, x_s)]ds \ge 0, \ t \in [v_{j},\tau_j],$$
proving that $y(t) \geq x(t)$ for all $t \in [\tau_j-h,\tau_j]$.
Now, for $t \in (\tau_j, s_j],$
\begin{eqnarray} \label{xtzt}
y(t)-x(t)=\int_{\tau_j}^t \exp\{-(t-s)\} (r(y(s-h))- f(s,
x_s))ds\leq 0,
\end{eqnarray}
since  
$f(s, x_s) \ge r(\mathcal{M}(x_s)) \geq r(\mathcal{M}(y_s)) =
r(y(s-h))$. 
Hence, by (\ref{xtzt}), $m_j=x(s_j)\ge y(s_j)\geq
F(M+\varepsilon)$ which proves that $m \geq F(M)$.
\end{proof}

To study the properties of $F$, we use its more explicit form
given below:
\begin{lemma}
\label{formula} Set $r^{-1}(u) = u/(a-u)$.  For $z \in I$,
$(F(z)-r(z))z \geq 0$ and
\end{lemma}
\begin{equation} \label{b6}
\theta = \int\limits_{r(z)}^{F(z)} \frac{ du}{r^{-1}(u) - r(z)}.
\end{equation}
\begin{proof} Let us consider $z > 0$, the case $z < 0$ being
similar. Consider the solution $y(t,z)$ of Eq. (\ref{ep}), recall
that $y(t,z) = z$ for $t \in [t_1-h, t_1],$ where $t_{1}=-\ln
(1-z/r(z))\in [-h,0]$. Next, $y(t,z)=r(z)(1-\exp{(-t)}),\, t\in
[t_1,t_1+h]$ so that $y(0,z) = 0$ and $y'(h,z)=-y(h,z)$. Therefore
$F(z)=y(t_*,z)$  at some point $t_*\in (t_1+h,h)$ where also
$y'(t_*,z)=0$.

Since $t_* \in [t_1+h, h]$, by the variation of  constants formula
we have
\begin{equation}
\label{mmm}
y(t_*,z)=F(z)=e^{-(t_*-h)}\left[y(t_1+h,z)e^{t_1}+\int\limits_{t_1}^{t_*-h}
e^v r(y(v,z))dv\right].
\end{equation}
On the other hand, $y'(t_*,z) = 0 =- y(t_*,z) + r(y(t_* - h,z))$,
so that $F(z)=y(t_*,z) = r(y(t_* - h,z)) \geq r(z)$ and $
r^{-1}(F(z))=y(t_*- h,z) = r(z)[1-\exp\{-(t_*- h)\}]. $ Thus
\begin{eqnarray*} \label{b1}
t_*-h =   \ln (r(z)/[r(z)-r^{-1}(F(z))]).
\end{eqnarray*}
Let now $y(v,z)=w$ (so that $\exp(v)=r(z)/(r(z)-w)$), then
\begin{eqnarray*}
 & &\int\limits_{t_1}^{t_*-h}e^vr(y(v,z))dv =
\int\limits_{z}^{r^{-1}(F(z))} r(w)d\frac{r(z)}{r(z)-w}  =\\[3mm]
& &= r(z)\left[\frac{r(z)}{z-r(z)}-\frac{F(z)}{r^{-1}(F(z))-r(z)}+
\int\limits_{z}^{r^{-1}(F(z))}\frac{dr(w)}{w-r(z)}\right].
\end{eqnarray*}
Putting now the last expression and the values of $t_1$, $t_*-h$
in (\ref{mmm}), we get (\ref{b6}).
\end{proof}

Finally, we state an important technical lemma whose proof can be
found in Appendix, Lemmas \ref{leform1}, \ref{plyus}.
\begin{lemma} \label{leform1F} Assume that
$(a, \theta) \in \mathcal{D}.$ Then $ F(z) <  \mathcal{R}(r(z))$
if $z \in ((a\beta -1)^{-1}, 0),$ and $ F(z) >  \mathcal{R}(r(z))$
 if $ z \in (0, a(\theta-1)/\theta - 1), $ where
$\mathcal{R}$ is defined in Subsection \ref{nodo}.
\end{lemma}

We will also consider  $\mathcal{F}: (a_*, +\infty)\to \R$ defined
by $\mathcal{F}(x) = F(x/(a-x))$. It can be easily seen that
$\mathcal{F}(r(z)) = F(z)$ for all $z \in I$.

\subsection{One-dimensional map $\mathbf{F_1: [0,+\infty) \to} (a,0]$}
By definition,  for $z\geq 0,$ $F_1 (z) = \min_{t \in [0,h]} y(t)$,
where
$y(t,z)$ satisfies (\ref{ep}) and has the initial value $y(s,z) =
(1 - e^{-s})r(z), \ s\in [-h, 0]$. We will need the following
\begin{lemma} \label{Func1}
Let $x(t)$ be a solution of (\ref{1pr}) and set $M =
\limsup_{t \to \infty}x(t),$ $m = \liminf_{t \to \infty}x(t)$. If
$(a,\theta) \in \mathcal{D}$, then $m \ge F_1 (M).$
\end{lemma}
\begin{proof}
Take
$\varepsilon, s_j, t_j, m_j, M_j, \tau_j$
as in the first two paragraphs of the proof of Lemma \ref{44}. Then,  for $t
\in [\tau_j-h, \tau_j],$ we have
\begin{eqnarray*}
x(t) =  \int_{\tau_j}^{t} e^{-(t-u)} f(u, x_u ) du \le
\int_{\tau_j}^{t} e^{-(t-u)} r(M+\varepsilon) du =
y(t-\tau_j,M+\varepsilon). \label{fu22}
\end{eqnarray*}
Thus, if $u \in [s_j-h, s_j]$, then $\mathcal{M}(x_u(s)) \leq
\mathcal{M}(y_u(s,M+\varepsilon))$ so that
\begin{eqnarray*} \label{fu33}
f(u, x_u) \ge   r( \mathcal{M}(y_u(s,M+\varepsilon)))=  r(r(M)(1 -
e^{-(u-h-\tau_j)})).
\end{eqnarray*}
This implies that
\begin{eqnarray*} \label{fu44}
 m_j &=& x(s_j) \ge \int_{\tau_j}^{s_j} e^{-(s_j-u)}r(r(M)(1 -
e^{-(u-h-\tau_j)}))du\\
&  = &\int_{\tau_j}^{s_j} e^{-(s_j-u)}r(y(u-h,M+\varepsilon))du =
y (s_j-\tau_j, M+\varepsilon ) \ge F_1 (M+\varepsilon).
\end{eqnarray*}
Since $\varepsilon > 0$ and $m_j \to m$ are arbitrary, the lemma
is proved.
\end{proof}
\begin{lemma} \label{form2} Set $r_1(z) =r(r(z)(1-e^h))$. For $z
> 0$ we have that $F_1(z)>a$ and
\begin{equation} \label{bb6}
\frac{r_1(z)\theta}{r(z)} = \int_{r_1(z)}^{F_1(z)} \frac{
du}{r^{-1}(u) - r(z)}.
\end{equation}
\end{lemma}
\begin{proof}
Take $t_*\in (0,h)$ such that
\begin{equation}
\label{efe} F_1(z)= y(t_*,z)= \int_{0}^{t_*}e^{-(t_* - u)}
r(r(z)(1 - e^{-(u-h)})) du.
\end{equation}
Since $y'(h) > 0$,
we have that $y'(t_*)=0$ and
therefore $F_1(z)=y(t_*) = -y'(t_*) + r(y(t_* - h))=r(y(t_* - h))>a.$
This implies that $ r^{-1}(F_1(z))=y(t_* - h) =
r(z)(1-\exp\{-(t_*- h)\}), $ from where
\begin{equation} \label{bb1}
t_* - h = \ln (r (z)/[r(z)-r^{-1}(F_1(z))]).
\end{equation}
Now, using (\ref{bb1}) and setting $\xi = r(z)(1 - e^{-(u-h)})$ in
(\ref{efe}), we obtain
\begin{eqnarray*}
 F_1(z)
&=& -(r(z) - r^{-1}(F_1(z))) \int_{r(z)(1 - e^h)}^{r^{-1}(F_1(z))}
r(\xi) d\frac{1}{\xi - r(z)}   \\[3mm]
& =&  -(r(z) - r^{-1}(F_1(z))) \Biggl( \frac{r(\xi)}{\xi - r(z)}
\Biggl|_{r(z)(1 - e^h)}^{r^{-1}(F_1(z))} + \int_{r(z)(1 -
e^h)}^{r^{-1}(F_1(z))} \frac{d r(\xi)}{\xi - r(z)} \Biggl).
\end{eqnarray*}
Simplifying this relation, we obtain (\ref{bb6}).
\end{proof}

We conclude this section by stating two lemmas which compare
$F_1$ and the associated function $\mathcal{F}_1(r)
\stackrel{def}= F_1(r/(a-r))$ with rational functions. The proofs
of these statements are based on rather careful estimations of
identity (\ref{bb6}) and are given in Appendix, Lemmas
\ref{leform2}, \ref{FuncRr2},  \ref{lele2} (it should be noted
that $\mathcal{R}$ approximates extremely well $\mathcal{F}_1$ so
that a very meticulous analysis of (\ref{bb6}) is needed).
\begin{lemma} \label{leform11F}  If $(a,\theta) \in \mathcal{D}^*$
and $z \ge a(\theta-1)/\theta -1$ then
$F_1(z) > \mathcal{R}(r(z)). $
\end{lemma}
\begin{lemma} \label{leform22F}
If $(a,\theta) \in \mathcal{S}$ and $z > 0$, then
\begin{equation} \label{F1}
\mathcal{F}_1(r(z)) >  \frac{1 + \ln\theta  - \theta}{2 +\ln\theta
- \theta}\frac{ar(z)}{1+ r(z)\frac{1 + \ln\theta - \theta}{1  -
\theta}} = \mathcal{R}_2(r(z)).
\end{equation}
Furthermore, $\mathcal{R}_2(a)> -1$ and $r(\mathcal{R}_2(a)) <
1/\beta.$
\end{lemma}
\subsection{Proof of Theorem \ref{main}}
 Let $x:[\alpha-h, \infty) \to \R$ be a solution of Eq.
(\ref{1pr}) and set $M = \limsup_{t \to \infty}x(t), \ m =
\liminf_{t \to \infty}x(t)$. We will reach a contradiction if we
assume that $m<0<M$ (note that the cases $M \leq 0$ and $m \geq 0$
were already considered in Lemma \ref{L1}).

First suppose that $(a,\theta) \in \mathcal{S}$. By Lemmas
\ref{Func1} and \ref{leform22F}, we obtain that
\begin{eqnarray} \label{pron}
m \ \ge \ F_1(M) = \mathcal{F}_1(r(M)) \ > \ \mathcal{R}_2(r(M))
> -1.
\end{eqnarray}
Take now an arbitrary $z \geq 0$.  Since $r(z) \in (a,0]$ and
$\mathcal{R}_2(z)$ is increasing on $(a,0]$, we get
$r(\mathcal{R}_2(r(z))) < 1/\beta$ due to Lemma \ref{leform22F}.
Therefore, the rational function $\lambda \stackrel{def} =
\mathcal{R}\circ r\circ \mathcal{R}_2\circ r: [0, \infty) \to [0,
\infty)$ is well defined. By Lemmas \ref{44} and \ref{leform1F},
we obtain
$$M \le \mathcal{F}(r(m)) < \mathcal{R}(r(\mathcal{R}_2(r(M)))) =
\lambda(M).$$  On the other hand, due to the inequality
$\lambda'(0) = \gamma(a,\theta) < 1$ (see Lemma \ref{albetand}),
we obtain that $\lambda(z) < z$ for all $z > 0$, a contradiction.

Let now $(a,\theta) \in \mathcal{D}^*$ and define the rational
function $R:[0, +\infty) \to (-\infty,0]$ as  $R =
\mathcal{R}\circ r$. We note that (\ref{global stability
condition}) implies $R'(0) = a\alpha(a,\theta) \in (-1,0)$. Next,
\begin{eqnarray} \label{uslo1}
m > R(M) > \mathcal{R}(a)  > -1.
\end{eqnarray}
Indeed,  if $M \le a(\theta-1)/\theta - 1,$ then Lemmas \ref{44}
and \ref{leform1F} imply that $ m \ge F(M) > \mathcal{R}(r(M)) =
R(M)$. If  $M \ge a(\theta-1)/\theta - 1,$ then Lemmas
\ref{Func1} and \ref{leform11F} give that $ m \ge F_1(M)
 > \mathcal{R}(r(M)) = R(M). $ The last inequality in (\ref{uslo1}) follows from
Lemma~\ref{albetand}.
Finally, applying Lemmas \ref{44} and \ref{leform1F}, and
using (\ref{uslo1}) and the inequality $R\circ R(x) < x , \ x >
0,$ which holds since $(R\circ R)'(0) = (R'(0))^2 < 1$, we obtain
that $$M \le F(m)  < \mathcal{R}(r(m)) < \mathcal{R}(r(R(M))) =
R(R(M)) < M,$$
 a  contradiction.

To prove the second part of Theorem \ref{main} take  $a < 0$ and
$h > 0$  such that (\ref{global stability condition}) is not satisfied.
Then, by Theorem 2.9 from \cite{ilt} there is a continuous functional $f$ satisfying (\ref{ycond}) and
such that the equilibrium $x(t)=0$ in (\ref{1pr}) is not locally
asymptotically stable.

\section{Some estimations of the global attractor for (\ref{LW3})}\label{Pr2}
To complete the proof of Theorem \ref{son}, we need to estimate
the bounds of the global attractor to (\ref{LW3}). We start by
stating a result from \cite{gyt}
\begin{lemma}{\rm \cite{gyt}} \label{luck} Let $q > 1$. Then there exist
finite positive limits $$M = \limsup_{t \to \infty}x(t)\, ,\, m =
\liminf_{t \to \infty}x(t)$$ for every nonnegative solution
$x(t)\not\equiv 0$ of (\ref{LW3}). Moreover, $[m,M] \subseteq
g([m,M])$ and $[m,M] \subseteq g_1([m,M]),$ where $g_1 = \theta
\ln q + (1-\theta)g,\ \theta = \exp(-\tau).$
\end{lemma}

Since the global stability of Eq. (\ref{LW3}) for  $\ln q \in
(0,2]$ was already proved in \cite{gyt}, we can suppose that $\ln
q
> 2$. In this case the minimal root $x_1$  of equation $g(x_1)= \ln q \ $
belongs to the interval $(0,1)$. Note that $x =1 < \ln q$ is the
point of absolute maximum for $g$ and $g_1 = (1-\theta)g + \theta
\ln q,$ so that  $g(1) > \ln q$ and $g_1(1) > \ln q$. We will use
the information about the values of $g$ and $g_1$ at $x_1 < 1 <
\ln q$ in the subsequent analysis.

Now, let us consider an arbitrary solution $x(t)$ of (\ref{LW3})
and its bounds $m, M$ defined in Lemma \ref{luck}. It is clear
that if we prove the existence of $m_*= m_*(q)$ such that $m \geq
m_*(q)
> x_1$ and $
m_*(q)$ does not depend on $x(t)$, then the change of variables $y
= x - \ln q$ transforms Eq. (\ref{LW3}) into an equation
satisfying {\bf (W)} within the domain of attraction, and
therefore Theorem \ref{main} can be applied.

Since $[m, M] \subseteq g([m, M])$, we obtain immediately that
either $m = M =\ln q$ or $m < \ln q < M.$ In the first case the
theorem is proved, so we will consider the second
possibility. Next, since $z < g(z)$ for $z \in (0, \ln q)$, we
have that $g(m) > m$ and
$$[m, M] \subseteq g([m, M]) \subseteq [\min\{g(m),g(M)\},g(1)]
= [g(M),g(1)].$$ Hence, $[m, M] \subseteq g([g(M),g(1)]) \subseteq
[\min\{g^2(1),g^2(M)\},g(1)]$. On the other hand, since $g(M) <
\ln q$ we get analogously that $g^2(M)> g(M)$. Next, since $g$ is
decreasing on $[1, +\infty)$ and $g(1) \geq M$ we find that
$g^2(1) \leq g(M)$. Thus $g^2(1) \leq g(M) < g^2(M)$ so that
$\min\{g^2(1),g^2(M)\}= g^2(1)$ and $[m, M] \subseteq
[g^2(1),g(1)]$. Therefore $m \ge g(g(1))$. Since the inequality $m
\ge g_1(g_1(1))$ can be proved analogously, the proof of theorem
will be completed if we establish that $m_*(q) =
\max\{g^2(1),g_1^2(1)\} > x_1$. We have

\noindent {\bf (i)} \underline{$g^2(1) > x_1$  for all $\ln q \in
[2,2.833157].$} This is an obvious fact if $g^2(1) \geq 1$, so that we only
need to consider the case $x_1, g^2(1)\in (0,1)$. Since $g$ is
increasing on $(0,1)$, the inequality  $g^2(1) > x_1$ is
equivalent to $g^3(1) > g(x_1) = \ln q$  in this case. Finally, a
direct computation shows that
$$g^3(1) - \ln q = q^3  e^{-1-q/e} \exp (-q^2 e^{-1-q/e}) - \ln q>0$$
 whenever $\ln q \in [2, 2.833157].$

\noindent {\bf (ii)} \underline{ $g_1^2(1) > x_1$  for all $\ln q
> 2.5$.} First, let us note that $x_1 \le \ln q + y_1,$ where
$$y_1 = \left(2-\ln q- \sqrt{(\ln q)^2 + 4\ln q - 4}\right)/2
$$ is the negative root of $\tilde{g}(y) = (y +
\ln q)(1 - y + y^2 / 2) - \ln q$. Indeed, with $x = y + \ln q$ and
$y \in (y_1,0)$, we have that $g(x)- \ln q = qxe^{-x} - \ln q = (y
+ \ln q)e^{-y} - \ln q \ge \tilde{g}(y) > 0.$

Since $g_1^2(1) \ge g_1(+\infty)= \theta \ln q,$ to finish the
proof of  {\bf (ii)}, it suffices to show that $\theta \ln q \ge
\ln q + y_1.$ Taking into account (\ref{globstab}) and using the
inequality $\ln(1+x) \ge x / (1+x)$, we obtain that
\begin{eqnarray*}(\theta-1)\ln
q -y_1&=&(\theta-1)(1+c) - y_1\geq (1+c)(-1 +c\ln ((c^2 +c)/(c^2 +
1)) - y_1\\ &\ge& -2 - \frac{2-\ln q- \sqrt{(\ln q)^2 + 4\ln q -
4}}{2}
\ge 0, \quad {\rm for} \quad \ln q \ge 5/2. \end{eqnarray*}

\section{Appendix}
\subsection{Preliminary estimations}
$\,$
\begin{lemma} \label{albet} For all $(a,\theta) \in \mathcal{D}$,
we have that $\alpha(a,\theta) > 0, \ \beta(a,\theta)> 0$, and
\end{lemma}
$$T(a,\theta)\stackrel{def}= \ (a^2 - a) \beta(a,\theta)(1 - \theta) +
\alpha(a,\theta) - (1-\theta) \ge 0.
$$
\begin{proof} Since $a(\theta -1) > 1$ for all
$(a,\theta) \in \mathcal{D}$ and
\begin{equation}\begin{array}{ll}
 1+x < \exp(x) < 1 + x + x^2/2 & \; {\rm for} \ x < 0, \\
\noalign{\smallskip}
\exp(x) > 1 + x + x^2/2+ x^3/6 & \; {\rm for} \ x > 0,
\end{array}
\label{expo}
\end{equation}
we have $ \alpha(a,\theta) = (1-a)\exp(\theta/a)+ a >
(1-a)(1+\theta/a) +a = 1 -\theta +\theta/a > 0.$ Analogously,
$\beta(a,\theta) > 0 $ for
all $(a,\theta) \in \mathcal{D}$ because of the following chain of
relations 
\begin{eqnarray*}
&- & a \alpha(a,\theta) \beta(a,\theta)e^{-\theta/a} =
  a^2e^{-\theta/a}- (1-a)^2e^{\theta/a}+ (1-2a +2\theta
(a-1))
\\[3mm]
& \geq&  a^2(1- \frac{\theta}{a}+\frac{\theta^2}{2a^2}-
\frac{\theta^3}{6a^3}) - (1-a)^2(1+
\frac{\theta}{a}+\frac{\theta^2}{2a^2}) + (1-2a +2\theta (a-1))  \\[3mm]
& =&  \frac{\theta}{6a^2}(-3\theta - a(\theta^2-6\theta +6)) >
\frac{\theta(2\theta a - 2a - \theta)}{2a^2}
> 0.
\end{eqnarray*}
To prove that $T(a,\theta) \ge 0$ for all $(a,\theta) \in
\mathcal{D}$, we replace $\alpha, \beta$ with their values in
$\alpha T$:
\begin{equation}
\begin{array}{ll}
\quad \alpha(a,\theta)T(a,\theta) =&\!\!\! a(2a-a^2-a\theta +
a^2\theta-1+\theta)
\\[2mm]
\, &\!\!\!+ 2\theta(a-1)(-2a+a\theta+1-\theta)\lambda
-(1-a)^2(-a+a\theta-\theta)\lambda^2.  \label{dop1}
\end{array}\end{equation}
It should be noticed that $-a+a\theta-\theta = a(\theta-1) - \theta
> 1-\theta>0$. Similarly, $2\theta(a-1)(-2a+a\theta+1-\theta) < 0$
so that $T(a,\theta) \ge 0$ if
$$\theta(-a+a\theta-\theta)(4a^4)^{-1}[-\theta(-\theta+a\theta-2a)^2
+ 4(\theta-1)a^3] > 0, 
$$
where the last expression was obtained from (\ref{dop1}) by
replacing $\lambda = \exp(\theta/a)$ by  $1 + \theta/ a +
\theta^2/2 a^2
> \exp(\theta /a)$.  Taking into account that $a(\theta-1)
> 1$ and $\theta < 1$ for $(a,\theta) \in \mathcal{D}$, we end the
proof of this lemma by noting that $ 4(\theta-1)a^3
-\theta(\theta-a\theta+2a)^2   \ge
 4a^2-(\theta-a\theta+2a)^2
 = (-4a - \theta + a\theta)\theta (1-a) > 0. $
\end{proof}
\begin{lemma} \label{albeta}
For all $(a,\theta) \in \mathcal{D}$ one has
\end{lemma}
\begin{equation}\label{alberta}
\frac{a \alpha(a,\theta)}{1- a\beta(a,\theta)} > -1.
\end{equation}
\begin{proof} It follows directly from the definitions of
$\alpha(a,\theta)$ and $\mathcal{D}$ that $a \alpha(a,\theta) >
-1$ for all $(a,\theta) \in \mathcal{D}$. Now, (\ref{alberta})
follows from the fact that $a\beta(a,\theta)< 0$ if $(a,\theta)
\in \mathcal{D}$.
\end{proof}
\begin{lemma}\label{dom} If $(a,\theta) \in \mathcal S,$ then
$\gamma(a,\theta) < 1.$
\end{lemma}
\begin{proof}
Notice that $(a,\theta) \in \mathcal S$ implies that $ \theta \in
[0.8,1)$ and $\theta>\Pi_3(-1/a)$ (or, equivalently,
$a>\pi_3(\theta)\stackrel{def}=-1/\Pi_3^{-1}(\theta)=(108 +
25\theta)(95\theta-95)^{-1}$). 
Here $\Pi_3^{-1}$ is the inverse function of $\Pi_3$.
Next we prove the inequality
\begin{equation}\label{ml}
\frac{1-\theta+ \ln \theta}{2-\theta+ \ln \theta} +
\frac{(\theta-1)^2(2-\theta)}{2} >0, \quad \theta \in [0.8,1),
\end{equation}
which is equivalent to the relation
$$
\Xi(q) \stackrel{def}= q(q-2)(q^2+1)+(2+q^2-q^3)\ln(q+1) > 0, \quad q \in
[-0.2,0),
$$
with $\theta = q+1$ (note that $1-q +\ln(q+1) > 0$ for $q \in
[-0.2,0)$). To do that, we will need the following approximation
of $ \ln(1+q)$ when  $q \in [-0.2,0)$:
\begin{eqnarray} \label{R303}
\ln(1+q) > q-0.5q^2+0.4q^3, \ \ q \in [-0.2,0).
\end{eqnarray}
(Indeed,  function $y(x)= \ln(1+x) -(x-0.5x^2+0.4x^3)$ has
exactly one critical point $x= -1/6$ on $[-0.2,0)$, and $y(-0.2) =
0.00005644\dots > 0$,  $y(0)=0$). Inequality (\ref{R303})
implies that
$
\Xi(q) \geq -0.1q^3(5q+2-9q^2+4q^3).
$
Now, since $(5q+2-9q^2+4q^3) > 0$ for all $q \in [-0.2,0)$,
we have that $\Xi(q) > 0$ and thus (\ref{ml}) is proved.

Next, due to (\ref{expo}) and (\ref{ml}),
$$
\gamma(a,\theta)
\le a^3(1 - \theta + \frac \theta a - \frac{\theta^2}{2 a} +
\frac{\theta^2}{2a^2})\frac{(\theta-1)^2(\theta-2)}{2} = w(a,\theta),
$$
so that $\gamma(\pi_3(\theta),\theta))\leq w(\pi_3(\theta),\theta)).$
Now, $w(\pi_3(\theta),\theta))$ is a fifth degree polynomial and an
elementary analysis shows that $w(\pi_3(\theta),\theta))<1$
 for all $\theta\in [0.8, 1).$
Since
\begin{eqnarray*}
& & \partial w(a,\theta)/\partial a =0.25(2 - \theta)(\theta-1)^2[
a(\theta-1)(2a+2\theta) + a( 4a(\theta-1)- 2\theta)] < 0,
\end{eqnarray*}
we conclude that $\gamma(a,\theta) \le w(a,\theta) < 1$ for $a >
\pi_3(\theta).$
\end{proof}
\begin{lemma}\label{Jcal}
Let $(a,\theta) \in \mathcal{D}$ and $r > -1/4$. Set
$\mathcal{J}(r) = \mathcal{I}(N(r))$, where $\mathcal{I}(N) =N
\coth(\nu N/2), \ \nu = -\theta/a, \ N(r) = \sqrt{1 + 4r}.$ Then
$\mathcal{J}'(0)
> 0$ and
\begin{equation}
\label{pFcal} \mathcal{J}(r) \le \mathcal{J}(0) + \mathcal{J}'(0)
r = \frac{1 + \lambda}{1 - \lambda} + \Biggl(2\frac{1 + \lambda}{1
- \lambda} + \frac{4\theta\lambda}{a(1 - \lambda)^2}\Biggl)r.
\end{equation}
\end{lemma}
\begin{proof} Set $k(N)= e^{2\nu N} - 2\nu N e^{\nu N} - 1 > 0$; then
$k(0) = 0$ and $k'(N) = 2\nu e^{\nu N}[e^{\nu N} - 1 - \nu N] > 0$
for all $N >0$. Hence $ \mathcal{I}\,^{\prime}(N)= k(N)/(e^{\nu N}
- 1)^2  > 0$ and $\mathcal{J}'(0) = \mathcal{I}\,'(1) N'(0)
> 0.$

Next, since $dN/dr = 2(1+4r)^{-1/2} = 2/ N, \quad d^2 N/dr^2 = -4
/ N^3,$ we obtain that
\begin{eqnarray*}
\mathcal{J}''(r) &=& \frac{\partial^2 \mathcal{I}(N(r))}{\partial
r^2} = \frac{\partial^2 \mathcal{I}(N)}{\partial N^2}
\Biggl(\frac{\partial N(r)}{\partial r}\Biggl)^2 + \frac{\partial
\mathcal{I}(N)}{\partial N}\Biggl(\frac{\partial^2 N(r)}{\partial
r^2}\Biggl)\\
 &=& \frac{4[- e^{3\nu N} + e^{2\nu N}(2\nu^2 N^2 - 2 \nu N +1) +
e^{\nu N}(2\nu^2 N^2 + 2 \nu N +1) - 1]}{N^3(e^{\nu N} - 1)^3}\\
 &=&\frac{\sum_{j=0}^{+\infty}p_j (\nu N)^j}{N^3(e^{\nu N} - 1)^3} <
 0, \quad {\rm since\ }  (\nu N) > 0, \ p_j = 0, \ j = 0,\dots , 5, {\rm \
 and}
\end{eqnarray*}
$$p_j = \frac{4}{(j-2)!}\Bigl(\frac{-3^j + 2^j + 1}{j(j-1)} +
\frac{-2^j + 2}{j-1} + 2^{j-1} + 2\Bigl) < 0,  \ j \ge 6.$$ Thus
$\mathcal{J}(r) \le \mathcal{J}(0) + \mathcal{J}'(0) r$ and
(\ref{pFcal}) is proved.
\end{proof}

\subsection{Properties of function $\mathbf{F}$}
To study the properties of functions $F: I \to \Bbb R$ and
$\mathcal{F}: J= (a_*, +\infty) \to \Bbb R,$ defined in subsection
\ref{subF}, it will be more convenient to use the integral
representation (\ref{b6}) instead of the original definition of
$F$. It should be noted that  conditions $xF(x) <0, \
(F(x)-r(x))x >0, \ x\in I\setminus \{0\}$ define $F$ in a unique
way: moreover $F$ and $\mathcal{F}$ are continuous and smooth at
$0$ with $\mathcal{F}'(0)= \alpha(a,\theta)$, $\mathcal{F}''(0)=
2\alpha(a,\theta)\beta(a,\theta)$. We have taken into consideration these
facts to define the rational functions $R$
and $\mathcal{R}$ (see Subsection \ref{subF}); however, since we do not use
anywhere these characteristics of $F$, their proof is omitted here.
\begin{lemma} \label{leform1} Assume that $r \in [a_*,0]$. Then
\end{lemma}
\begin{equation}
\mathcal{F}(r) > \alpha\, r/(1 - \beta r) = \mathcal{R}(r).
\end{equation}
\begin{proof}
1. First, suppose that $4r + 1 > 0$ and $r \not= 0$. Since $0 >
\mathcal{F}(r)
> r$, we have, for every $z \in [r, \mathcal{F}(r)]$ and $a < 0$,
\begin{equation}\label{ta1}
r^{-1}(z) = z/(a - z)  = (z/a) (1 + (z/a) + (z/a)^2 + \cdots)> (z/a)
+ (z/a)^2.
\end{equation}
Hence
\begin{equation} \label{ta2}
\theta = \int_r^{\mathcal{F}(r)} \frac{dz}{z(a-z)^{-1}-r} <
\int\limits_{r}^{\mathcal{F}(r)} \frac{dz}{\frac za + (\frac za)^2
- r} = a \int\limits_{r/a}^{\mathcal{F}(r)/a} \frac{du}{u + u^2 -
r}. \end{equation} Now, since for $r < 0$ the roots $\alpha_1 =
(-1 - \sqrt{1 + 4r})/2, \quad \alpha_2 = (-1 + \sqrt{1+ 4r})/2$ of
the equation $ u^2 + u - r =0 $ are negative, we obtain
\begin{eqnarray}
 a \int\limits_{r/a}^{\mathcal{F}(r)/a} \frac{du}{u + u^2 - r}
\label{urav12}
 =  -\frac{a}{\sqrt{1 + 4r}} \ln \Bigl( \frac{\mathcal{F}(r) -
a\alpha_1} {\mathcal{F}(r) - a\alpha_2} \ \frac {r - a\alpha_2}{r
- a\alpha_1} \Bigl) > \theta.
\end{eqnarray}
 The last inequality implies that
\begin{equation} \label{urav2}
\frac{\mathcal{F}(r) - a\alpha_1}{\mathcal{F}(r) - a\alpha_2} \
\frac{r - a\alpha_2}{r - a\alpha_1} > \exp (\frac{\theta\sqrt{1 +
4r}}{-a})
 \stackrel{def}= \omega(r) \ge 1.
 \end{equation}
Taking into account that $\mathcal{F}(r) - a\alpha_2 < 0$ and $r -
a\alpha_1 < 0$, and replacing $\alpha_1, \alpha_2$ in (\ref{urav2})
by their values,  we obtain
\begin{eqnarray*}
\mathcal{F}(r) > \frac{r(2a^2 - a + a\frac{\omega(r) + 1}{\omega(r) -
1}\sqrt{1 + 4r})} {2r + a + a\frac{\omega(r) + 1}{\omega(r) - 1}\sqrt{1
+ 4r}}.
\end{eqnarray*}
Next, since $\mathcal{J}(r) = \frac{\omega(r) + 1}{\omega(r) - 1}\sqrt{1
+ 4r}$, we can apply  Lemma \ref{Jcal} to see that
\begin{equation}\begin{array}{ll}
 \mathcal{F}(r) &\!\!\!> \displaystyle\frac{r(2a - 1 +
\mathcal{J}(r))}{2r/a + 1 +
\mathcal{J}(r)} \ \geq \ \displaystyle\frac{r(2a - 1 +
\mathcal{J}(0) + \mathcal{J}'(0)r)}{2r/a + 1 + \mathcal{J}(0)
+ \mathcal{J}'(0)r}  \nonumber \\\noalign{\medskip}
\, &\!\!\!= \displaystyle\frac{r(\lambda + a(1 - \lambda)) + \frac12
\mathcal{J}'(0) (1 - \lambda) r^2} {1 + (\frac{1 - \lambda}{a} + \frac12
\mathcal{J}'(0) (1 - \lambda)) r} \stackrel{def}=\mathcal{L}(r).
\label{guk}
\end{array}\end{equation}
Now,  $\mathcal{L}(r)=(a_1r  + a_2 r^2)/(1 + a_3 r)$, with
$a_1 > 0, a_2 > 0$. Moreover, since $0 < \mathcal{J}(r) \le
\mathcal{J}(0) +
\mathcal{J}'(0)r,$ all denominators in (\ref{guk}) are positive so
that  $1 + a_3 r > 0$. Next,
\begin{equation} \label{p20}\qquad
 a_1 a_3 - a_2=(\lambda + a(1 - \lambda)) \Biggl(\frac{1 -
\lambda}{a} + \frac12 \mathcal{J}'(0) (1 - \lambda)\Biggl) -
\frac12 \mathcal{J}'(0) (1 - \lambda)\le 0.
\end{equation}
Indeed, the last inequality is equivalent to the obvious relation
\begin{eqnarray*}
\frac{\lambda + a(1 - \lambda)}{a} < 0 \le \frac12 (1 - \lambda)(1
- a) \mathcal{J}'(0)
\end{eqnarray*}
(notice that $\alpha= \lambda + a(1 - \lambda) > 0,$ while, by
Lemma \ref{Jcal}, $\mathcal{J}'(0) > 0$).

Finally, since  the inequality
$ (a_1 r + a_2 r^2)(1 + a_3 r)^{-1} \ge a_1 r(1 + (a_3 - a_2/a_1)
r)^{-1}$ holds for $r < 0, \ a_1
> 0,\ a_2 > 0, \ 1 + a_3 r > 0,\ a_1 a_3 \le a_2 $,  we obtain
\begin{eqnarray*}
\mathcal{F}(r) > \frac{a_1 r}{1 + (a_3 - \frac{a_2}{a_1}) r} =
 \frac{r(\lambda + a(1 - \lambda))}{1 + (\frac{1 -
\lambda}{a} + \mathcal{J}'(0)\frac{1 - \lambda}{2}- \frac{
\mathcal{J}'(0) (1 - \lambda)}{2(\lambda + a(1 - \lambda))})r} =
\mathcal{R}(r).
\end{eqnarray*}
Hence the statement of the lemma is proved for $r \in (-1/4,0)$.
As an important consequence of the first part of proof, we get the
following relation  $$ \lim\limits_{r \to -1/4} \frac{r(2a^2 - a +
a\frac{\omega(r) + 1}{\omega(r) - 1}\sqrt{1 + 4r})} {2r + a +
a\frac{\omega(r) + 1}{\omega(r) - 1}\sqrt{1 + 4r}} =
 \frac{a^2 (1 - \theta) + \theta a/2}{\theta(2a-1) - 4a^2} \ge
 \mathcal{R}(-1/4),$$ which will be used in the next stage of proof.

2. The case $r = -1/4.$ 

From (\ref{ta2}), evaluated at $r=-1/4$,
we get $\theta <(-2a^2) / (2\mathcal{F}(-1/4)+a) + (4a^2)/(2a-1),$
so that
$$ \mathcal{F}(-1/4) > \frac{a^2 (1 - \theta) + \theta a/2}{\theta
(2a-1) - 4a^2} \ge \mathcal{R}(-1/4). $$

3. Assume now that $4r + 1 < 0.$ We have
\begin{equation}
 \qquad\quad
a \int\limits_{r/a}^{\mathcal{F}(r)/a} \frac{du}{u + u^2 - r}
=
\frac{2a}{\sqrt{-4r-1}}\Biggl(\arctan\frac{2\mathcal{F}(r)+a}{a\sqrt{-4r-1}}
- \arctan\frac{2r+a}{a\sqrt{-4r-1}}\Biggl). \label{le1}
\end{equation}
By (\ref{ta2}) and (\ref{le1}), we obtain
\begin{eqnarray*}
2\mathcal{F}(r) + a >
a\sqrt{-4r-1}\frac{\frac{2r+a}{a\sqrt{-4r-1}} + \tan\frac{\theta
\sqrt{-4r-1}}{2a}}{1- \frac{2r+a}{a\sqrt{-4r-1}} \tan\frac{\theta
\sqrt{-4r-1}}{2a}}.
\end{eqnarray*}
 Now, since  $\tan x \le x + x^3 / 3$ for $x \in (-\pi/2, 0)$
and $a < 0$,  we obtain
\begin{eqnarray} \label{gra}
\mathcal{F}(r) > r \frac{a^2(1-\theta) + \theta a/2 + \theta ^3(-r
- \frac 14)(\frac{1}{2a}
  - 1)(3)^{-1}}{a^2 - \theta(r + \frac a2) - \theta^3
(-r - \frac 14)(\frac a2 + r)(3a^2)^{-1}} = G(r).
\end{eqnarray}
Therefore it will be sufficient to establish that $G(r) \ge
\mathcal{R}(r)$  for $r < -1/4$. First, note that by the second
part of the proof
\begin{eqnarray*}
G(-1/4) = \frac{a^2 (1-\theta) + \theta a/2}{\theta(2a-1) - 4a^2}
\ge \mathcal{R}(-1/4).
\end{eqnarray*}

Let us consider now the function $H(r)= G(r) - \mathcal{R}(r)$ for
$r \le 0.$ Since $G(r) = G_1(r) / G_2(r)$, where $G_j$ are
polynomials in $r$ of second degree, $H(r)$ can be written  as
\begin{equation} \label{h1}
H(r) = \frac{G_1(r)(1-\beta r) - \alpha r G_2(r)}{G_2(r)(1-\beta
r)} = \frac{H_1(r)}{H_2(r)},
\end{equation}
so that $H$ is a quotient of two polynomials of  third degree with
$H_2(r) >0$ for $r \le 0.$  We get $\lim_{r \to -\infty}
G(r) = a^2(1 - 1/(2a))
> 0,$ and therefore $H(-\infty)= \lim_{r \to -\infty} H(r) > 0.$
Furthermore,  $H(0) = 0$ and
\begin{eqnarray*}
H'(0) = \frac{1-\theta(1 - \frac{1}{2a}) +
\frac{\theta^3}{12a^2}(1- \frac{1}{2a})} {1 - \frac{\theta}{2a} +
\frac{\theta ^3}{24a^3}} - (a + e^{\frac{\theta}{a}} (1 - a)) =
\frac{\sum_{k=5}^{+\infty}p_k\theta^k} {1 - \frac{\theta}{2a} +
\frac{\theta ^3}{24a^3}} > 0,
\end{eqnarray*}
since the denominator of the last fraction is positive and
$p_{2m+1}> 0, \ p_{2m} < 0, \ p_{2m+1} + p_{2m+2} > 0$ for $m \ge
2$. Here we use the formula
$$p_k = \frac{a-1}{a^k k!}(1 - \frac k2 + \frac{k(k-1)(k-2)}{24}),
\ k \ge 5. $$
 Finally, since $H(-1/4)
= G(-1/4) - \mathcal{R}(-1/4) \ge 0$, there exists at least one
zero of $H(r)$ in the interval $[-1/4, 0)$.  $H_1(r)$ is a
polynomial of third degree in $r$ and therefore it cannot have more than three
zeros. Hence, since $H(-\infty)
> 0$ and $H(-1/4) \ge 0,$ we obtain that $H(r) \ge 0$ if $r <-1/4.$
\end{proof}
\begin{lemma} \label{plyus} If $(a,\theta) \in \mathcal{D}$, then
$\mathcal{F}(r) < \mathcal{R}(r)$ for all $r \in (0, 1/\beta)$.
\end{lemma}
\begin{proof} By definition of $\mathcal{F}$, we have that
$z > 0$ if $r > 0$ and $z \in [\mathcal{F}(r), r]$. We begin the
proof  by assuming that $r \in (0, a^2-a).$ Since $z
> 0, a < 0$, we find that $z(a - z)^{-1} < z/a +z^2/a^2$.
Therefore, (\ref{ta2}) holds under our present conditions. Now, since
$r>0$, the roots of  equation $u^2 + u - r = 0$ are $ \alpha_1 = (-1 -
\sqrt{1 + 4r})/2<0, \quad
\alpha_2 = (-1 +
\sqrt{1 + 4r})/2>0.$
Next, since $\mathcal{F}(r) -
a\alpha_1  < r - a\alpha_1 < 0$  for all $r \in (0, a^2-a)$, we
obtain that the relations (\ref{urav12}),
(\ref{urav2}) hold in the new situation and therefore
\begin{eqnarray*}
\mathcal{F}(r) < \frac{a^2 \alpha_1 \alpha_2 (\omega(r) - 1) + ar
(\alpha_1 - \omega(r) \alpha_2)} {- a\alpha_2 + r + \omega(r)(
a\alpha_1 -r)}=\frac{r(2a^2 - a + a\frac{\omega(r) + 1}{\omega(r) -
1}\sqrt{1 + 4r})} {2r + a + a\frac{\omega(r) + 1}{\omega(r) - 1}\sqrt{1
+ 4r}},
\end{eqnarray*}
where the denominator is positive for every $r \in (0, a^2-a)$.
Now, recall that $\mathcal{J}(r) = \frac{\omega(r) + 1}{\omega(r) -
1}\sqrt{1 + 4r}$ ; applying Lemma \ref{Jcal}, we obtain
$\mathcal{J}(r) \le \mathcal{J}(0) + \mathcal{J}'(0) r$. Next,
since  for all $r \in (0,a^2 -a)$ we have that $2r/a + 1 +
\mathcal{J}(r)= -2(a(\omega(r)-1))^{-1}(- a\alpha_2 + r + \omega(r)(
a\alpha_1 -r))
> 0,$ and  the function $p(x)= (r(2a - 1 + x))(2r/a + 1 +
x)^{-1}$ is increasing in $x$, we get $\mathcal{F}(r)<\mathcal{L}(r)$
(compare with (\ref{guk})).
Now, $ (a_1 r + a_2 r^2)(1 + a_3 r)^{-1} \le a_1 r(1 + (a_3 -
a_2/a_1) r)^{-1}$ if  $a_1 a_3 - a_2 \leq 0,\ r > 0,$ $  \ a_1 >
0,\ a_2 > 0$ . Therefore, by (\ref{p20}), $\mathcal{L}(r) \leq
\mathcal{R}(r)$.

Now we assume that $r \ge a^2-a.$ Taking into account that
$z(a-z)^{-1} < 0$ for $z > 0,$ we obtain the  inequality $$ \theta =
\int_{r}^{\mathcal{F}(r)} \frac{dz}{z(a-z)^{-1} - r} <
\int_{r}^{\mathcal{F}(r)} \frac{dz}{-r} = \frac{\mathcal{F}(r)-r}{-r},$$
so that $\mathcal{F}(r) < r(1-\theta)$. Finally, the inequality
$r(1-\theta) \le \mathcal{R}(r) = \alpha r/(1-\beta r)$ is
equivalent to $ r \ge (1- \theta - \alpha)/((1 - \theta) \beta)$,
which holds for all $r \ge a^2-a $ due to the relation $a^2-a \geq (1-
\theta - \alpha)/((1 - \theta) \beta)$, established in Lemma
\ref{albet}.
\end{proof}

\subsection{Properties of function $\mathbf{F_1}$ in the domain $\mathcal{D}^*$}
Suppose now that $(a,\theta) \in \mathcal{D}^*.$ We study some
properties of function $F_1$ and  the associated function
$\mathcal{F}_1: (a, 0) \to \R$ defined as $\mathcal{F}_1(r(z)) =
F_1(z),$ which, by Lemma \ref{form2}, satisfies
$$
\frac{r_1(r)\theta}{r} = \int\limits_{r_1(r)}^{\mathcal{F}_1(r)}
\frac{ dz}{r^{-1} (z) - r}\, ,\, \mbox{ where } r_1(r) = \frac{ar(\theta-1)}{\theta +r(\theta-1)}.
$$
\begin{lemma} \label{leform2}  Assume that $(a,\theta) \in
\mathcal{D}^*$ and that the inequalities
$ a < r \leq a_* = a + \theta/(1-\theta)$ hold. Then  $\mathcal{F}_1(r) >
\mathcal{R}(r). $
\end{lemma}
\begin{proof}
Since $r_1(r)<\mathcal{F}_1(r)<0$ and
$a_* < -1$ for $(a,\theta) \in \mathcal{D}^*,$  using (\ref{ta1})
we get $$\frac{r_1(r)\theta}{r} <
\int\limits_{r_1(r)}^{\mathcal{F}_1(r)} \frac{dz}{\frac za +
(\frac za)^2 - r} = a \int\limits_{r_1(r)/a}^{\mathcal{F}_1(r)/a}
\frac{du}{u + u^2 - r}. $$ The last integral can be transformed as
it was done in  (\ref{le1}) to obtain
\begin{eqnarray*}
\frac{ r_1(r)\theta}{r} <
\frac{2a}{\sqrt{-4r-1}}\Biggl(\arctan\frac{2\mathcal{F}_1(r)+a}
{a\sqrt{-4r-1}} - \arctan\frac{2r_1 + a}{a\sqrt{-4r-1}}\Biggl).
\end{eqnarray*}
Therefore
$$ \varsigma _1 \stackrel{def}=
\arctan\frac{2\mathcal{F}_1(r)+a}{a\sqrt{-4r-1}} < \arctan\frac{2r_1
+ a}{a\sqrt{-4r-1}} + \frac{\theta r_1\sqrt{-4r-1}}{2ar} \stackrel{def}=
\varsigma_2+ \varsigma_3,$$ and since $\varsigma_1, \varsigma_2
\in (0,\pi/2), \varsigma_3 < 0$, we obtain
\begin{eqnarray*}
2\mathcal{F}_1(r) + a > a\sqrt{-4r-1}\frac{\frac{2r_1 +
a}{a\sqrt{-4r-1}} + \tan\frac{\theta r_1 \sqrt{-4r-1}}{2ar}}{1-
\frac{2r_1 + a}{a\sqrt{-4r-1}} \tan\frac{\theta r_1
\sqrt{-4r-1}}{2ar}}.
\end{eqnarray*}
Since $\tan x < x + x^3 / 3$  for $x \in (-\pi/2, 0)$, we have
\begin{eqnarray} \label{vopr}
\mathcal{F}_1(r) > \frac{A_1(P)r + A_2(P)r^2}{B_0(P) + B_1(P)r +
B_2(P)r^2}
 = G_1 (r, P, a, \theta),
\end{eqnarray}
where
\begin{eqnarray*}
& & A_1(P) = (1 - \theta)P + \frac{\theta}{2a}P^2 +
\frac{\theta^3}{24a^3}(2a - P)P^3, \quad A_2(P) = \frac{\theta^3}{6a^3}(2a - P)P^3, \\
& & B_0(P) = 1 - \frac{\theta P}{2a} + \frac{\theta^3 P^3}{24a^3},
\quad B_1(P) = -\frac{\theta P^2}{a^2} + \frac{\theta^3 P^3}{6a^3}
+ \frac{\theta^3 P^4}{12 a^4}, \quad B_2(P) = \frac{\theta^3
P^4}{3a^4}, \\
& & P = P(r,a,\theta) = r_1 / r = \frac{a(\theta-1)}{\theta +
r(\theta - 1)}.
\end{eqnarray*}
After substitution of the value of $P$ into (\ref{vopr}), we get
\begin{eqnarray*}
\mathcal{F}_1(r) > G_1(r, P(r,a,\theta),a,\theta) \stackrel{def}=
\mathcal{G}_1 (r,a,\theta) = \frac{r M(r,a,\theta)}{N(r,a,\theta)},
\end{eqnarray*}
 where\begin{eqnarray*}
 M(r,a,\theta) &=& 24(A_1(P) +
A_2(P)r)(\theta + r(\theta - 1))^4
\\  & =& -(\theta - 1)^2 a [13 \theta^3 - \theta^5 - 2 \theta^2
(\theta - 1)(\theta + 3)(3\theta - 8)r - \\ & &-
 4 \theta (2\theta^2 - 15)(\theta - 1)^2 r^2 + 24 (\theta -
1)^3 r^3], \\[1mm]
 N(r,a,\theta) &=& 24(B_0(P) +B_1(P)r + B_2(P)r^2)(\theta +
r(\theta - 1))^4 =  35 \theta^4 - 9 \theta^5 \\ & &+ \theta^7 -
3 \theta^6 + \theta^3 (\theta - 1 )(7\theta^3 - 17\theta^2-
47\theta + 153) r    + 12\theta^2 (\theta^3 - 2\theta^2  \\ & &-
6\theta + 19)(\theta - 1)^2 r^2  - 12\theta (3\theta - 11)(\theta
- 1)^3 r^3 + 24 (\theta - 1)^4 r^4.
\end{eqnarray*}
To prove our lemma, it suffices to check the inequality
$\mathcal{G}_1 (r,a,\theta) \geq \mathcal{R}(r)$ for $r \in [a,
a_*]$. First, considering $N(r,a,\theta)=N(r,\theta)$ as a
polynomial in $r$ of the form
$N(r,a,\theta)=\sum_{k=0}^{4}N_k(\theta) r^k$, we can check that
$(-1)^k N_k(\theta)>0$ for $\theta\in (0,1)$ and therefore,   for
all $ \theta\in (0,1)$ and $  r<0$,
\begin{equation}\label{nz}
N(r,a,\theta) = 24(B_0(P) +B_1(P)r + B_2(P)r^2)(\theta + r(\theta
- 1))^4>0.
\end{equation}
Since $N(r,a,\theta) > 0, \ 1 - \beta r > 0$ (recall that
$\beta(a,\theta) > 0$ in the domain $\mathcal{D}$), the inequality
$r M(r,a,\theta)/N(r,a,\theta) \geq \alpha r/(1-r\beta)$ is
equivalent to
\begin{eqnarray} \label{pol}
Q(r,a,\theta) \stackrel{def}= (1 - r\beta(a,\theta)) M(r,a,\theta) -
\alpha (a,\theta) N(r,a,\theta) \leq 0.
\end{eqnarray}
Now, an easy comparison  of  $\mathcal{G}_1(a_*,a,\theta) =
G_1(a_*,1,a,\theta)$ with $G(a_*)$ given in (\ref{gra}) shows that
the inequality (\ref{pol}) is fulfilled for $r = a_*.$ In next two
lemmas, we will prove that $\partial Q(r,a,\theta) /
\partial r > 0$ for all $r \in [a,a_*].$ Therefore, since
$Q(a_*,a,\theta) \leq 0,$ we obtain  $Q(r,a,\theta) \leq 0$
for $r \in [a,a_*],$ which proves that $\mathcal{F}_1(r) >
\mathcal{R}(r).$
\end{proof}

\begin{lemma}
\label{lele}$S(r,a,\theta) = \frac{\partial}{\partial
r}Q(r,a,\theta)
> 0$ at the point $r = a_*$.
\end{lemma}
\begin{proof}
Recall that we are interested in the case $r= a_*$, when $P=1$. By (\ref{pol}) and
the above definitions of $M(r,a,\theta), N(r,a,\theta)$,
$$ Q(r,a,\theta) = 24(\theta + r(\theta - 1))^4((A_1(P) +
A_2(P)r)(1 - \beta r) -
\alpha (B_0(P) + B_1(P)r + B_2(P) r^2)).
$$
Next, setting $P' = \partial P(r,a,\theta)/\partial r\mid_{r=a_*}
= -a^{-1}$, $A'_j = \partial A_j(a,\theta,P)/\partial P|_{P=1}, $
$B'_j =
\partial B_j(a,\theta,P)/\partial P|_{P=1}, \ A_j =A_j(a,\theta,1) ,
B_j=B_j(a,\theta,1)$, we obtain that
\begin{eqnarray}
& & \partial Q(r,a,\theta) / \partial r |_{r=a_*} = 24(
Q_1(r,a,\theta)+Q_2(r,a,\theta)),\label{Q12}
\end{eqnarray}
where
\begin{eqnarray*}
 Q_1 &=& 4 a^3(\theta-1)^4 ((A_1 +
A_2a_*)(1 - \beta a_*) - \alpha (B_0 + B_1 a_* + B_2 a_*^2))\\
& & + a^4(\theta - 1)^4((A'_1 + A'_2a_*)(1 - \beta a_*) -
\alpha(B'_0 + B'_1 a_* + B'_2 a_*^2))P';\\
 Q_2 &=&a^4 (\theta - 1)^4 (A_2 - \beta A_1 - \alpha B_1 - 2a_*
\beta A_2 - 2a_* \alpha B_2).
\end{eqnarray*}
Now, for the convenience of the reader, the following part of the
proof will be divided in several steps.

{\it \underline{Step (i): $Q_2(r,a,\theta) > 0$}}. Indeed,
consider the second degree polynomial
$$\chi_1(r) \stackrel{def}=
(A_1 + A_2r)(1 - \beta r) - \alpha (B_0 + B_1 r + B_2 r^2).$$
Notice that $ \chi_1(r)= \frac{H_1(r)}{r}, $ where $H_1$ is
defined in (\ref{h1}). This implies that the unique critical point
of $\chi_1$ belongs to $(-1/4, +\infty)$ and that
$\chi_1(+\infty) = -\infty$. Hence $\chi'_1(r) > 0$ for all $r <
-1/4$  so that $Q_2(\theta, a, r) =
  a^4(\theta - 1)^4 \chi'_1(r) > 0$. 

{\it  \underline{Step (ii):}} The following inequality holds:
\begin{equation}
(A'_1 + A'_2a_*)(B_0 + B_1a_* + B_2a_*^2) - (A_1 + A_2a_*)(B'_0 + B'_1a_* +
B'_2a_*^2) > 0.
\label{proverka1}
\end{equation}
Indeed,  the left-hand side of (\ref{proverka1}) can be
transformed into
\begin{equation}\begin{array}{l} 
  \displaystyle\frac{1}{576 a^6(1-\theta)^{3}}(-\theta^{6}(3\theta+1)^{3}
+ 12\theta^{6}(3\theta+1)^{2}(\theta-1)a  \\ \noalign{\smallskip}
 - 24\theta^4(\theta-1)(3\theta+1)(2\theta^3 - 2\theta^2 +
3\theta - 5)a^2 + 32\theta^3(2\theta^4 - 2\theta^3 + 18\theta^2
 \\
    - 39\theta - 9)(\theta -
1)^2a^3 - 48\theta^2(8\theta^3 - 41\theta^2 + 30\theta - 9)(\theta -
1)^2a^4  \\
 - 576\theta (\theta^2 - \theta +2)(\theta-1)^3 a^5 +
576(\theta-1)^4 a^6). \label{prov2}
\end{array}\end{equation}
Taking into account that $\eta \stackrel{def}= (\theta - 1)a > 1,$ the sum
of the first two terms in (\ref{prov2}) is positive:
$$ -\theta^{6}(3\theta+1)^{3} +
12\theta^{6}(3\theta+1)^{2}(\theta-1)a =  \theta^{6}(3\theta+1)^{2}(-(3\theta+1)
+ 12(\theta - 1)a) > 0.
$$
The other terms in (\ref{prov2}), can be written
as
\begin{eqnarray*}
& & a^2(- 24\theta^4(\theta-1)(3\theta+1)(2\theta^3 - 2\theta^2 +
3\theta - 5)\\ & & + 32\theta^3(2\theta^4 - 2\theta^3 + 18\theta^2
- 39\theta
- 9)(\theta -1) \eta \\
& & - 48\theta^2(8\theta^3 - 41\theta^2 + 30\theta - 9)\eta^2  - 576\theta
(\theta^2 - \theta +2)\eta^3 + 576 \eta^4) \stackrel{def}= a^2
\Upsilon(\theta, \eta).
\end{eqnarray*}
By  the Taylor formula,
\begin{eqnarray} \label{Taylor}
\Upsilon(\theta, \eta) = \Upsilon(\theta, 1) + (\eta-1)\partial
\Upsilon(\theta, 1)/\partial \eta  + 0.5 (\eta-1)^2
\partial^2\Upsilon(\theta, \eta_1)/\partial \eta^2,
\end{eqnarray}
where $\eta_1 \in [1, 2].$ It is easy to verify that
\begin{eqnarray*}
\Upsilon(\theta, 1) &=& -8(\theta-1)^2 [\theta^4
(18\theta^3-2\theta^2+27\theta-81) + (108\theta^3 -54\theta^2-72)]
> 0,\\
\partial \Upsilon(\theta, 1)/\partial \eta &=&
32(\theta-1)[(2\theta^6+18\theta^4+36)(\theta-1) -
45\theta^2(\theta-1)^2 - 36] > 0, \\
\partial^2 \Upsilon(\theta,\eta)/\partial \eta^2 &=&
3456\eta( 2\eta -\theta(\theta^2-\theta+2)) - 96\theta^2
(8\theta^3-41\theta^2+30\theta-9) > 0
\end{eqnarray*}
(here we use the inequality $ 8\theta^3-41\theta^2+30\theta-9 < 0,
\ \theta \in [0,1]$).

Finally, by (\ref{Taylor}), $\Upsilon(\theta, \eta) > 0$ for
$\theta \in (0,1), \eta \in (1,2).$ \vspace{3mm}

{\it \underline{Step (iii):}} We have
\begin{eqnarray} \label{prov3}
\varrho \stackrel{def}= (B'_0 + B'_1a_* + B'_2a_*^2)(B_0 + B_1a_* +
B_2a_*^2)^{-1} < 1.
\end{eqnarray}
Indeed, taking into account (\ref{nz}), the latter inequality is
equivalent to$$(\upsilon \stackrel{def}=) \qquad
\frac{\theta^3}{12a^3} + a_*(-\frac{\theta}{a^2}+
\frac{\theta^3}{3a^3} + \frac{\theta^3}{4a^4}) +
\frac{\theta^3}{a^4}a^2_* < 1.$$ Now, we know that $|a_*| \leq
|a|$ and $|a^{-1}| < 1 - \theta$ (so that $\theta/|a| < 1/4$).
Therefore $\upsilon < 1/4 + (1/3)(1/16) + 1/16 < 1$. \vspace{4mm}

{\it \underline{Step (iv): $Q_1(r,a,\theta) > 0.$}} First, using
(\ref{nz}) and (\ref{proverka1}), we obtain that
\begin{eqnarray*}
(A'_1 + A'_2a_*)(1 - \beta a_*) - \alpha(B'_0 + B'_1 a_* + B'_2
a_*^2) \\
\ge \varrho ((A_1 + A_2a_*)(1 - \beta a_*) - \alpha (B_0 + B_1 a_*
+ B_2 a_*^2)).
\end{eqnarray*}
 Next, using inequality (\ref{pol}) which was
proved at $r= a_*$, we find that
$$(A_1 + A_2a_*)(1 -
\beta a_*) - \alpha (B_0 + B_1 a_* + B_2 a_*^2) =
\frac{Q(a_*,a,\theta)}{24(\theta + a_*(\theta-1))^4}< 0.$$
Therefore,
\begin{eqnarray*}Q_1(r,a,\theta) \geq a^3(\theta - 1)^4((A_1 + A_2a_*)(1
- \beta a_*) - \alpha (B_0 + B_1 a_* + B_2 a_*^2))(4-\varrho) > 0.
\end{eqnarray*}
{\it \underline{Step (v):}} Recalling (\ref{Q12}) and Steps {\it
(i),(iv)}, we finish the proof of the lemma.
\end{proof}
\begin{lemma} \label{leleka}
$S(r,a,\theta) > 0$ for $r \in [a, a_*].$
\end{lemma}
\begin{proof} Differentiating function $Q$ given by (\ref{pol}), we obtain
\begin{eqnarray*}
S(r,a,\theta) &=& \sum\limits_{i=0}\limits^{3} S_i (\theta,a) r^i = 96
(\theta-1)^4 (\beta a(\theta-1) - \alpha)r^3  \\[1mm] & & +
(-12(\theta-1)^4 a\theta (2\theta^2 - 15)\beta + 36 \theta (\theta - 1)^3
(3\theta - 11)\alpha - 72 (\theta-1)^5 a ) r^2
\\[1mm]
& & + ( -4 a \theta^2 (\theta-1)^3 (\theta+3)(3\theta-8)\beta -
24\theta^2 (\theta-1)^2 (\theta^3 - 2\theta^2 - 6\theta +
19)\alpha  \\[1mm]
& & + 8 a \theta (\theta-1)^4 (2\theta^2 - 15))r  \\[1mm]
& & +  a(\theta-1)^2 \theta^3 (13 - \theta^2)\beta - \theta^3
(\theta-1)(7\theta^3 - 17\theta^2 - 47\theta + 153)\alpha  \\[1mm]
& & + 2 a \theta^2 (\theta - 1)^3 (\theta+3)(3\theta - 8).
\end{eqnarray*}
Now, inequalities $S_i (a,\theta)r^i < 0, \ i = 3, 2, 1, 0,$
are equivalent to
$a\theta (\theta-1)\beta > T_i (a,\theta),$ where
\begin{eqnarray*}
& & T_3(a,\theta) = \theta\alpha, \quad  T_2(a,\theta) = \frac{-6
a (\theta-1)^2 +
3\theta(3\theta - 11)\alpha}{2\theta^2 - 15}, \label{S2} \\[2mm]
& & T_1(a,\theta) = \frac{2a(\theta-1)^2 (2\theta^2 - 15) -
6\theta (\theta^3 - 2\theta^2 - 6\theta +
19)\alpha}{3\theta^2 + \theta - 24}, \label{S1} \\[2mm]
& & T_0(a,\theta) =  \frac{2 a (\theta-1)^2 (3\theta^2 + \theta - 24) - \theta (7
\theta^3 - 17 \theta^2 - 47 \theta + 153)\alpha }{\theta^2 - 13}.
\label{S0}
\end{eqnarray*}

Next, for $(a,\theta)  \in  \mathcal{D}^*, $ the following
inequalities hold
\begin{eqnarray}
T_3(a,\theta) > T_2(a,\theta), \label{RA3} \\[2mm]
T_2(a,\theta) > T_1(a,\theta), \label{RA1} \\[2mm]
T_1(a,\theta) > T_0(a,\theta) \label{RA2}.
\end{eqnarray}
Indeed,
taking into account that $\alpha  = (1-a)\exp(\theta/a)+a$,
 inequality
(\ref{RA1}) is equivalent to
\begin{equation}\begin{array}{l}
 4\theta^6 - 8\theta^5 - 41\theta^4 + 108\theta^3 + 99\theta^2
- 312\theta - 306 +  \\[1mm]
 +3\theta (4\theta^5- 8\theta^4- 45\theta^3  + 106\theta^2 +
97\theta -306)e^{\theta / a} (1-a)a^{-1} <0.\label{r21a}
\end{array}\end{equation}
Since $3\theta (4\theta^5- 8\theta^4- 45\theta^3  + 106\theta^2 +
97\theta -306)<0$,
it is sufficient to prove (\ref{r21a}) for the maximum value in
$a$ of the function $\frac{1-a}{-a} e^{\frac{\theta}{a}}.$ The
derivative of this function is equal to $-e^{\theta / a}(a\theta -
a - \theta)/ a^3, $ and it is positive if $a < \theta /(\theta -
1).$ Hence, it is sufficient to verify (\ref{r21a}) at $a =
\pi_1(\theta) = -1/\Pi_1^{-1}(\theta) = (1 + \sqrt{1 +
4\theta(1-\theta)})/(2(\theta-1))$ if $\theta \in (0, 0.8],$ and at
$a = \pi_3(\theta) =-1/\Pi_3^{-1}(\theta)= (133 +
25(\theta-1))/(95(\theta-1))$ if $\theta \in [0.8, 1).$ 

 Using
(\ref{expo}) and replacing the value $a = \pi_3(\theta)$ in
(\ref{r21a}), we get the following  expression:
\begin{eqnarray*}
& & \frac{q^2}{2(133+25q)^3}[-40522972 + 220135634 q - 410248779
q^2 + 204446752 q^3 + 279016108 q^4  \\
& &  + 23396520 q^5 - 209505145 q^6
- 30804850 q^7 + 35072100 q^8 + 7581000 q^9],
\end{eqnarray*}
which is negative for $\theta = q+1 \in [0.8, 1).$ Direct
computations show that (\ref{r21a}) holds if $a = \pi_1(\theta)$
and $\theta \in[0, 0.8].$

Analogously, inequality (\ref{RA3}) is equivalent to
\begin{eqnarray} \label{r32a}
& & -(6 -3\theta^2 +6\theta + 2\theta^3) - \theta (18 - 9\theta +
2\theta^2) e^{\theta / a} (1-a)a^{-1}  < 0.
\end{eqnarray}
Using (\ref{expo}) and substituting the value $a = \pi_3(\theta)$ in
(\ref{r32a}), we get the expression
\begin{eqnarray*}
& & \frac{q^2}{2(133+25q)^3}[-4465209 - 971090 q - 12743680 q^2
\\
& & - 5731130 q^3 - 2242475 q^4 + 1103900 q^5 - 1263500 q^6],
\end{eqnarray*}
which is negative for $\theta = q+1 \in [0.8, 1).$ Direct
computations show again that (\ref{r32a}) is satisfied for $a =
\pi_1(\theta)$ and $\theta \in[0, 0.8].$

Finally, (\ref{RA2}) is equal to
\begin{equation}\begin{array}{l}
 (\theta^6 - 16\theta^5 + 2\theta^4 + 226\theta^3 + 63\theta^2
- 570\theta - 762)+     \\ \noalign{\smallskip}+ \theta
(15\theta^5 - 32\theta^4 - 212\theta^3 + 550\theta^2 +
813\theta-2190)e^{\theta / a} (1-a)a^{-1} < 0. \label{r10a}
\end{array}\end{equation}
Next, employing (\ref{expo}) and using the value $a = \pi_3(\theta)$ in
(\ref{r10a}), we get the expression
\begin{eqnarray*}
& & \frac{q^2}{2(133+25q)^3}[-14595952 + 471367808 q - 1571124744
q^2  - 18258802 q^3 + 723267159 q^4 \\
& & + 399356020 q^5 - 311046000 q^6
-73063100 q^7 + 42576625 q^8 + 9476250 q^9],
\end{eqnarray*}
which is negative for $\theta = q+1 \in [0.8, 1).$ Direct
computations also show in this case that (\ref{r10a}) holds if $a =
\pi_1(\theta)$ and $\theta \in[0, 0.8].$

To finish the proof of this lemma, we take 
an arbitrary $r \in [a, a_*]$ (so that $r = a_* k, k \geq 1$) and
 write  function $S(r,a,\theta)$ in the form
\begin{eqnarray*}
S(r,a,\theta) = \sum_{i=0}^{3} S_i (a,\theta) a_*^i k^i = k^3 (S_3
a_*^3 + \frac 1k S_2 a_*^2 + \frac 1{k^2} S_1 a_* + \frac 1{k^3}
S_0).
\end{eqnarray*}
First, note that $S_3 a_*^3
> 0$. Indeed,
if $S_3 a_*^3 \leq 0,$ then, in view of (\ref{RA3})-(\ref{RA2}),
$S_i a_*^i \leq 0$ for $i=0,1,2$, and therefore $S(a_*,a,\theta)
\leq 0$, contradicting  Lemma \ref{lele}. Next, the conclusion
of Lemma \ref{leleka} is obvious if $S_i a_*^i \geq  0$ for all
$i=0,1,2$. Finally, if $S_i a_*^i \leq 0$ and $S_{i+1} a_*^{i+1} >
0$ for some $i$, then using the above representation for
$S(r,a,\theta)$ and relations (\ref{RA1})-(\ref{RA2}), it is easy
to see that $S(r,a,\theta) \geq S(a_*,a,\theta)>0$ for $r \in [a,
a_*].$
 \end{proof}

\subsection{Properties of function $\mathbf{F_1}$ in the domain $\mathcal{S}$}
$\,$

\begin{lemma} \label{FuncRr2} If $ r\in (a,0 )$ and $h \le 1$, then
\begin{equation}
\label{F} \mathcal{F}_1(r) >
\frac{1 - h - e^{-h}}{2 - h - e^{-h}}\frac{ar }{1+ r\frac{1 - h -
e^{-h}}{1  - e^{-h}}} =  \mathcal{R}_2(r).
\end{equation}
\end{lemma}
\begin{proof} Take $z > 0$ and consider the point $t_* \in (0,h)$ defined
in (\ref{efe});  by Lemma \ref{form2}, $F_1(z) > a$. Since
$r(r(z)(1 - e^{-(s-h)})) < 0$ for all $s \in (0,h)$, it follows
from (\ref{efe}) that
\begin{equation}\begin{array}{l}
 \mathcal{F}_1(r(z)) =  F_1(z) >
 e^{-(t^*-h)} \displaystyle\int_0^h e^{s-h}r(r(z)(1 - e^{-(s-h)}))ds=
\\
[6mm]
  =\displaystyle\frac{r(z) - r^{-1}(F_1(z))}{r(z)} 
\displaystyle\int_{-h}^0 e^u r(r(z)(1 - e^{-u}))du
=   \phi (r(z)) - r^{-1}(F_1(z)) \psi(r(z)), \end{array} \label{f1}
\end{equation}
where
$$ \psi(x) =
\phi(x)/x, \ \phi(x) = \int_{-h}^0 e^u r(x(1 - e^{-u}))du.$$
Applying Jensen's inequality \cite[p.~110]{roy} to the last
integral, we obtain that
\begin{equation}\begin{array}{ll}
\phi(x) &\!\!\!\! =\displaystyle\int_{-h}^0 (1 - e^{-h}) r(x(1 -
e^{-u}))d\left(e^{u}/(1 - e^{-h})\right)\\
\, &\!\!\!\!\ge (1 - e^{-h}) r\Biggl( \displaystyle\frac{\int_{-h}^0
x(e^u - 1) du}{1 - e^{-h}}\Biggl)
 = \displaystyle\frac{ax(1 - h - e^{-h})}{1 + x\frac{1 -
h - e^{-h}}{1 - e^{-h}}}\stackrel{def} = x\mathcal{H}(x) .
\label{jen}
\end{array}\end{equation}
Denote $\psi=\psi(r)$, $\phi=\phi(r)$,
$\mathcal{F}_1=\mathcal{F}_1(r)$. Now, for $r < 0$, (\ref{f1})
implies that $\mathcal{F}_1  > \phi - (\mathcal{F}_1 \psi)/(a -
\mathcal{F}_1 )$. Since $a - \mathcal{F}_1  < 0,$ we conclude that
\begin{eqnarray}
\mathcal{F}_1^2 - \mathcal{F}_1 (\phi + \psi + a) + a\phi > 0.
\label{vvv}
\end{eqnarray} Next we prove that, under our assumptions,
\begin{eqnarray}
(\psi + \phi + a)^2 - 4a \phi > (\psi + \phi - a - 2 \psi_0)^2
\geq 0, \label{www}
\end{eqnarray}
where $\psi_0 = a(1 - h - e^{-h})$. Indeed, (\ref{www}) amounts to
\begin{eqnarray*}
\psi(\psi_0r + a + \psi_0) > \psi_0(a + \psi_0). 
\end{eqnarray*}
Since $\psi_0r + a + \psi_0 < 0$, the latter inequality is
equivalent to
\begin{eqnarray*}
\psi < \frac{\psi_0(a + \psi_0)}{\psi_0r + a + \psi_0} =
 \frac{a(1 - h - e^{-h})}{1 + r\frac{1 - h - e^{-h}}{2 - h -
 e^{-h}}}\,  \stackrel{def} = \mathcal{G}(r),
\end{eqnarray*}
which holds because for $a < 0, r < 0, h \leq 1$ we have
$\mathcal{H}(r) <   \mathcal{G}(r),$ and since, by (\ref{jen}), $
\psi(r) \le \mathcal{H}(r)$. Now, the inequalities $a\phi(r(z)) >
0$, (\ref{www}) and the continuous dependence of $\phi(r), \psi(r),
\mathcal{F}_1(r)$ on $r \in (a,0)$ imply that the quadratic
polynomial $y(x)= x^2 - x(\phi + \psi + a) + a\phi$ has two roots
$x_1 =x_1(r) < x_2 =x_2(r)$ with the same sign and that this sign
is the same for all $r  \in (a,0)$.
Similarly, by (\ref{vvv}), we have that  either $\mathcal{F}_1(r) <
x_1(r)$ or $\mathcal{F}_1(r) > x_2(r)$ for all $r \in (a,0)$.
Since
 $\mathcal{F}_1(0^-) = 0 > x_1(0^-) = \psi_0 +a$, we
conclude that $x_1(r), x_2(r)$ are negative for all $r  \in
(a,0)$, and $\mathcal{F}_1(r)
> x_2(r)$. In  other words,
\begin{equation}\begin{array}{ll}
 \mathcal{F}_1 &\!\!\!\!> \frac12 (\psi + \phi + a + \sqrt{(\psi +
\phi + a)^2 -
4a \phi}) \\ \noalign{\medskip}
\, &\!\!\!\!=   \displaystyle\frac{2a\phi}{\psi + \phi + a - \sqrt{(\psi
+ \phi + a)^2 - 4a \phi}} \ge \displaystyle\frac{2a\phi}{2(a + \psi_0)},
\label{piip}
\end{array}\end{equation} 
where
the last inequality is due to the following consequence of
(\ref{www}):
$$ \sqrt{(\psi + \phi + a)^2 - 4a \phi} \ge -a + \phi + \psi -
2\psi_0.$$
Finally, combining (\ref{jen}) and (\ref{piip}), we
obtain (\ref{F}).
\end{proof}
\begin{lemma} \label{lele2} Assume that 
$(a,\theta) \in \mathcal{S}$. Then
$r( \mathcal{R}_2(a)) < \beta^{-1}.$
\end{lemma}
\begin{proof} {\it \underline{Step (i)}:}
In the new variables
$q = \theta - 1, k = a(\theta -1)$,
the expression for $\mathcal{R}_2(a)$ takes the form
$$
\mathcal{R}_2(a) = \frac{-q +
\ln(q+1)}{1-q+\ln(q+1)}\frac{k^2}{q^2 - k(-q + \ln(q+1))}.$$
Next we prove that
\begin{eqnarray} \label{R301}
\mathcal{R}_2(a) \ge \frac{6k^2(q-1)}{3kq^2 - 4kq +12 +6k} \stackrel{def}=
\bar R_2(q,k)
\end{eqnarray}
for all $q \in [-0.2,0), k \in [1, 1.5].$ (Note that for $(a,\theta)\in
\mathcal{D}$, the inequalities $1\leq a(\theta-1)\leq 1.5$ hold).
Indeed, we have
\begin{eqnarray*}
\mathcal{R}_2(a) - \bar R_2 = \frac{-k^2
C(q,k)}{(1-q+L)(-q^2-kq+kL)(3kq^2-4kq+12+6k)},
\end{eqnarray*}
where $L = \ln(1+q),$ and
\begin{eqnarray*}
 C(q,k) &=& q(6q^3-12q^2+3kq^2+6q-8kq-12) \\  & &+(14kq+12-9kq^2+6q^2-6q^3)L +
6k(q-1)L^2.
\end{eqnarray*}
Next, the following inequalities hold in an obvious way for  $q
\in [-0.2,0)$ and $ k \in [1, 1.5]$:
\begin{eqnarray*} & & 1 - q +
\ln(1+q) \geq 1- q + q/(1+q) > 0,\\
& & 3kq^2 - 4kq + 12 +6k > 0,\\
& & -q^2 - kq + k\ln(1+q) < -q^2 < 0.
\end{eqnarray*}
Thus  $C(q,k)
> 0$ will imply that $\mathcal{R}_2(a)
> \bar R_2$. 
Now, in view of (\ref{R303}) and the obvious inequalities
$\min_{\{q \in [-0.2,0), k \in  [1, 1.5]\}}(14kq+12-9kq^2+6q^2-6q^3) \geq
7.26
> 0$ and $6k(q-1) < 0$, we obtain  that $C(q,k)>(q^3/50)  (-168k q^3 + 48
kq^4 - 120 q^3 +255
kq^2
+   270q^2 -150q-110kq-50k-60)\geq
-0.356176q^3 >0$.

{\it \underline{Step (ii)}:} Using the new variables, we obtain the
following expression for $\alpha$:
$$\alpha(a,\theta) = \alpha(k/q,1+q) = (1 -
k/q)\exp(q(q+1)/k) + k/q.$$ We will prove that $ \alpha > -q(24k^2
- 12k - 7q)/(24 k^2) \stackrel{def}= \bar\alpha
> 0.$  Indeed, since $ \exp(x) >
1+x+x^2/2+ x^3/6$ for all $x = q(1+q)/k < 0$, we get
$$
 \alpha(q,k) - \bar\alpha >  q^2(24k)^{-3}[4q^{4} + (12-4k)q^{3} +
12q^{2} + (-12k^2 + 12k + 4)q + k],
$$
where the right-hand side is positive for all $q \in [-0.2,0), k
\in [1, 1.5]$.

{\it \underline{Step (iii)}:} Set $E = \exp(\frac{q(q+1)}{k})$.   Here
we  prove that
\begin{eqnarray*}
 \alpha\beta &=& -k/q - (q/k)(-2q + 2k -
1)E + (k-q)^2(kq)^{-1}E^2  \nonumber \\[3mm]
& <&  \bar\beta_0 \stackrel{def}= q^2(q+1)[2k^2 - q(k+2k^2) +
q^2(9-11k+2k^2)](6k^4)^{-1}.
\end{eqnarray*}
Indeed, due to (\ref{expo}) and since  $\frac{-q}{k}(-2q + 2k - 1)
> 0,$ $ \frac{(k-q)^2}{kq} < 0,$  we obtain
\begin{eqnarray*}
\bar\beta_0 - \alpha\beta > (1/6)(q/k)^4(q+1)(-8q^2+10kq-16q+1)
> 0
\end{eqnarray*}
if $q \in [-0.2,0), \ k \in [1, 1.5]$.

{\it \underline{Step (iv)}:} First, note that $\bar R_2(q,k) > -1$ for all $q \in
[-0.2,0), k \in [1, 1.5]$ so that $r(\bar R_2(q,k))$ and
$r(\mathcal{R}_2(a))$ are well defined. Moreover, since $r$ is
strictly decreasing over $(-1,0)$, in virtue of (\ref{R301}) we
get
\begin{eqnarray} \label{R305}
0< r(\mathcal{R}_2(a)) < r(\bar R_2) = \frac{6k^3(q-1)}{q[3kq^2 +
q(6k^2 - 4k) -6k^2+6k+12]}.
\end{eqnarray}

{\it \underline{Step (v)}:} The above steps imply that 
$$r(\mathcal{R}_2(a))\beta = r(\mathcal{R}_2(a))\alpha\beta/\alpha<
r(\bar R_2) \bar\beta_0 /
\bar\alpha .$$ Hence, Lemma \ref{lele} will be proved if we show
that $r(\bar R_2) \bar\beta_0 / \bar\alpha  - 1 < 0$. We have
\begin{equation}\label{640}
\qquad r(\bar R_2) \frac{\bar\beta_0}{ \bar\alpha}  - 1 =
\frac{\sum_{i=0}^4 Z_i
q^i}{(3kq^2-4qk+12+6k-6k^2+6qk^2)(-12k-7q+24k^2)}\, ,
\end{equation}
where $ Z_0= 24k(6k^3-7k^2-9k+6),\ Z_1
=-144k^4+120k^3-114k^2+42k+84, \
  Z_2 =  -2k(36k^2+93k-94), \ Z_3 = 3k(16k^2+8k+7), \ Z_4 =-
24k(k-1)(2k-9).$

Now, in view of (\ref{R305}), we have that the denominator of the right-hand side of (\ref{640})
is positive for all $q \in [-0.2,0), k \in [1,
1.5]$.  Therefore it suffices to prove that
$\sum_{i=0}^4 Z_i q^i < 0$; we finish the proof 
by observing that, for $q \in [-0.2,0), \ k \in [1, 1.5],$
\begin{eqnarray*}
Z_0 + Z_1 q \leq Z_0 -0.2 Z_1 = 0.3(2k-3)(288k^3+112k^2-154k-5)-21.3< 0,\\
Z_2 q^2 =  -2kq^2(36k^2+93k-94)< 0, \\
Z_3+ Z_4 q \geq Z_3 - 0.2 Z_4  = 57.6k^3-28.8k^2+64.2k
> 0.
\end{eqnarray*}
\end{proof}

\section*{Acknowledgements} This research was
supported by FONDECYT (Chile), project 8990013. E. Liz was
supported in part by M.C.T. (Spain) and FEDER,  under project
BFM2001-3884-C02-02.
  V. Tkachenko was supported in part by
 F.F.D. of Ukraine, project 01.07/00109. The
authors are greatly indebted to an anonymous referee for his/her
valuable suggestions which helped them to improve the exposition
of the results.

\end{document}